\documentclass[11pt]{article}
\usepackage{amsmath, upgreek}
\usepackage{amssymb}
\usepackage{amscd}
\usepackage{xspace}
\usepackage{verbatim}
\usepackage{color}
\newcommand{\mysection}[1]{
\section{#1}\setcounter{equation}{0}}
\title{\bf Nonlinear elliptic equations with measure valued absorption potential}
\author{{\bf Nicolas Saintier\footnote{  Departamento de Matem\'atica, Facultad de Ciencias Exactas y Naturales, Universidad de Buenos Aires, Buenos Aires, Argentina. \newline 
Email: nsaintie@dm.uba.ar}}\qquad
 {\bf Laurent V\'eron\footnote{ D\'epartement de Math\'ematiques,  Universit\'e Fran\c{c}ois Rabelais, Tours, France.
 \newline Email: veronl@lmpt.univ-tours.fr }}\\[2mm]
}

\date{}
\begin{document}
 \maketitle


\newcommand{\txt}[1]{\;\text{ #1 }\;}
\newcommand{\tbf}{\textbf}
\newcommand{\tit}{\textit}
\newcommand{\tsc}{\textsc}
\newcommand{\trm}{\textrm}
\newcommand{\mbf}{\mathbf}
\newcommand{\mrm}{\mathrm}
\newcommand{\bsym}{\boldsymbol}
\newcommand{\scs}{\scriptstyle}
\newcommand{\sss}{\scriptscriptstyle}
\newcommand{\txts}{\textstyle}
\newcommand{\dsps}{\displaystyle}
\newcommand{\fnz}{\footnotesize}
\newcommand{\scz}{\scriptsize}
\newcommand{\be}{\begin{equation}}
\newcommand{\bel}[1]{\begin{equation}\label{#1}}
\newcommand{\ee}{\end{equation}}
\newcommand{\eqnl}[2]{\begin{equation}\label{#1}{#2}\end{equation}}
\newcommand{\barr}{\begin{eqnarray}}
\newcommand{\earr}{\end{eqnarray}}
\newcommand{\bars}{\begin{eqnarray*}}
\newcommand{\ears}{\end{eqnarray*}}
\newcommand{\nnu}{\nonumber \\}
\newtheorem{subn}{\name}
\renewcommand{\thesubn}{}
\newcommand{\bsn}[1]{\def\name{#1}\begin{subn}}
\newcommand{\esn}{\end{subn}}
\newtheorem{sub}{\name}[section]
\newcommand{\dn}[1]{\def\name{#1}}   
\newcommand{\bs}{\begin{sub}}
\newcommand{\es}{\end{sub}}
\newcommand{\bsl}[1]{\begin{sub}\label{#1}}
\newcommand{\bth}[1]{\def\name{Theorem}
\begin{sub}\label{t:#1}}
\newcommand{\blemma}[1]{\def\name{Lemma}
\begin{sub}\label{l:#1}}
\newcommand{\bcor}[1]{\def\name{Corollary}
\begin{sub}\label{c:#1}}
\newcommand{\bdef}[1]{\def\name{Definition}
\begin{sub}\label{d:#1}}
\newcommand{\bprop}[1]{\def\name{Proposition}
\begin{sub}\label{p:#1}}
\newcommand{\R}{\eqref}
\newcommand{\rth}[1]{Theorem~\ref{t:#1}}
\newcommand{\rlemma}[1]{Lemma~\ref{l:#1}}
\newcommand{\rcor}[1]{Corollary~\ref{c:#1}}
\newcommand{\rdef}[1]{Definition~\ref{d:#1}}
\newcommand{\rprop}[1]{Proposition~\ref{p:#1}}
\newcommand{\BA}{\begin{array}}
\newcommand{\EA}{\end{array}}
\newcommand{\BAN}{\renewcommand{\arraystretch}{1.2}
\setlength{\arraycolsep}{2pt}\begin{array}}
\newcommand{\BAV}[2]{\renewcommand{\arraystretch}{#1}
\setlength{\arraycolsep}{#2}\begin{array}}
\newcommand{\BSA}{\begin{subarray}}
\newcommand{\ESA}{\end{subarray}}
\newcommand{\BAL}{\begin{aligned}}
\newcommand{\EAL}{\end{aligned}}
\newcommand{\BALG}{\begin{alignat}}
\newcommand{\EALG}{\end{alignat}}
\newcommand{\BALGN}{\begin{alignat*}}
\newcommand{\EALGN}{\end{alignat*}}
\newcommand{\note}[1]{\textit{#1.}\hspace{2mm}}
\newcommand{\Proof}{\note{Proof}}
\newcommand{\qeda}{\hspace{10mm}\hfill $\square$}
\newcommand{\qed}{\\
${}$ \hfill $\square$}
\newcommand{\Remark}{\note{Remark}}
\newcommand{\modin}{$\,$\\[-4mm] \indent}
\newcommand{\forevery}{\quad \forall}
\newcommand{\set}[1]{\{#1\}}
\newcommand{\setdef}[2]{\{\,#1:\,#2\,\}}
\newcommand{\setm}[2]{\{\,#1\mid #2\,\}}
\newcommand{\mt}{\mapsto}
\newcommand{\lra}{\longrightarrow}
\newcommand{\lla}{\longleftarrow}
\newcommand{\llra}{\longleftrightarrow}
\newcommand{\Lra}{\Longrightarrow}
\newcommand{\Lla}{\Longleftarrow}
\newcommand{\Llra}{\Longleftrightarrow}
\newcommand{\warrow}{\rightharpoonup}
\newcommand{
\paran}[1]{\left (#1 \right )}
\newcommand{\sqbr}[1]{\left [#1 \right ]}
\newcommand{\curlybr}[1]{\left \{#1 \right \}}
\newcommand{\abs}[1]{\left |#1\right |}
\newcommand{\norm}[1]{\left \|#1\right \|}
\newcommand{
\paranb}[1]{\big (#1 \big )}
\newcommand{\lsqbrb}[1]{\big [#1 \big ]}
\newcommand{\lcurlybrb}[1]{\big \{#1 \big \}}
\newcommand{\absb}[1]{\big |#1\big |}
\newcommand{\normb}[1]{\big \|#1\big \|}
\newcommand{
\paranB}[1]{\Big (#1 \Big )}
\newcommand{\absB}[1]{\Big |#1\Big |}
\newcommand{\normB}[1]{\Big \|#1\Big \|}
\newcommand{\produal}[1]{\langle #1 \rangle}

\newcommand{\thkl}{\rule[-.5mm]{.3mm}{3mm}}
\newcommand{\thknorm}[1]{\thkl #1 \thkl\,}
\newcommand{\trinorm}[1]{|\!|\!| #1 |\!|\!|\,}
\newcommand{\bang}[1]{\langle #1 \rangle}
\def\angb<#1>{\langle #1 \rangle}
\newcommand{\vstrut}[1]{\rule{0mm}{#1}}
\newcommand{\rec}[1]{\frac{1}{#1}}
\newcommand{\opname}[1]{\mbox{\rm #1}\,}
\newcommand{\supp}{\opname{supp}}
\newcommand{\dist}{\opname{dist}}
\newcommand{\myfrac}[2]{{\displaystyle \frac{#1}{#2} }}
\newcommand{\myint}[2]{{\displaystyle \int_{#1}^{#2}}}
\newcommand{\mysum}[2]{{\displaystyle \sum_{#1}^{#2}}}
\newcommand {\dint}{{\displaystyle \myint\!\!\myint}}
\newcommand{\q}{\quad}
\newcommand{\qq}{\qquad}
\newcommand{\hsp}[1]{\hspace{#1mm}}
\newcommand{\vsp}[1]{\vspace{#1mm}}
\newcommand{\ity}{\infty}
\newcommand{\prt}{\partial}
\newcommand{\sms}{\setminus}
\newcommand{\ems}{\emptyset}
\newcommand{\ti}{\times}
\newcommand{\pr}{^\prime}
\newcommand{\ppr}{^{\prime\prime}}
\newcommand{\tl}{\tilde}
\newcommand{\sbs}{\subset}
\newcommand{\sbeq}{\subseteq}
\newcommand{\nind}{\noindent}
\newcommand{\ind}{\indent}
\newcommand{\ovl}{\overline}
\newcommand{\unl}{\underline}
\newcommand{\nin}{\not\in}
\newcommand{\pfrac}[2]{\genfrac{(}{)}{}{}{#1}{#2}}

\def\ga{\alpha}     \def\gb{\beta}       \def\gg{\gamma}
\def\gc{\chi}       \def\gd{\delta}      \def\ge{\epsilon}
\def\gth{\theta}                         \def\vge{\varepsilon}
\def\gf{\phi}       \def\vgf{\varphi}    \def\gh{\eta}
\def\gi{\iota}      \def\gk{\kappa}      \def\gl{\lambda}
\def\gm{\mu}        \def\gn{\nu}         \def\gp{\pi}
\def\vgp{\varpi}    \def\gr{\rho}        \def\vgr{\varrho}
\def\gs{\sigma}     \def\vgs{\varsigma}  \def\gt{\tau}
\def\gu{\upsilon}   \def\gv{\vartheta}   \def\gw{\omega}
\def\gx{\xi}        \def\gy{\psi}        \def\gz{\zeta}
\def\Gg{\Gamma}     \def\Gd{\Delta}      \def\Gf{\Phi}
\def\Gth{\Theta}
\def\Gl{\Lambda}    \def\Gs{\Sigma}      \def\Gp{\Pi}
\def\Gw{\Omega}     \def\Gx{\Xi}         \def\Gy{\Psi}

\def\CS{{\mathcal S}}   \def\CM^+{{\mathcal M}}   \def\CN{{\mathcal N}}
\def\CR{{\mathcal R}}   \def\CO{{\mathcal O}}   \def\CP{{\mathcal P}}
\def\CA{{\mathcal A}}   \def\CB{{\mathcal B}}   \def\CC{{\mathcal C}}
\def\CD{{\mathcal D}}   \def\CE{{\mathcal E}}   \def\CF{{\mathcal F}}
\def\CG{{\mathcal G}}   \def\CH{{\mathcal H}}   \def\CI{{\mathcal I}}
\def\CJ{{\mathcal J}}   \def\CK{{\mathcal K}}   \def\CL{{\mathcal L}}
\def\CT{{\mathcal T}}   \def\CU{{\mathcal U}}   \def\CV{{\mathcal V}}
\def\CZ{{\mathcal Z}}   \def\CX{{\mathcal X}}   \def\CY{{\mathcal Y}}
\def\CW{{\mathcal W}} \def\CQ{{\mathcal Q}}
\def\BBA {\mathbb A}   \def\BBb {\mathbb B}    \def\BBC {\mathbb C}
\def\BBD {\mathbb D}   \def\BBE {\mathbb E}    \def\BBF {\mathbb F}
\def\BBG {\mathbb G}   \def\BBH {\mathbb H}    \def\BBI {\mathbb I}
\def\BBJ {\mathbb J}   \def\BBK {\mathbb K}    \def\BBL {\mathbb L}
\def\BBM {\mathbb M}   \def\BBN {\mathbb N}    \def\BBO {\mathbb O}
\def\BBP {\mathbb P}   \def\BBR {\mathbb R}    \def\BBS {\mathbb S}
\def\BBT {\mathbb T}   \def\BBU {\mathbb U}    \def\BBV {\mathbb V}
\def\BBW {\mathbb W}   \def\BBX {\mathbb X}    \def\BBY {\mathbb Y}
\def\BBZ {\mathbb Z}

\def\GTA {\mathfrak A}   \def\GTB {\mathfrak B}    \def\GTC {\mathfrak C}
\def\GTD {\mathfrak D}   \def\GTE {\mathfrak E}    \def\GTF {\mathfrak F}
\def\GTG {\mathfrak G}   \def\GTH {\mathfrak H}    \def\GTI {\mathfrak I}
\def\GTJ {\mathfrak J}   \def\GTK {\mathfrak K}    \def\GTL {\mathfrak L}
\def\GTM {\mathfrak M}   \def\GTN {\mathfrak N}    \def\GTO {\mathfrak O}
\def\GTP {\mathfrak P}   \def\GTR {\mathfrak R}    \def\GTS {\mathfrak S}
\def\GTT {\mathfrak T}   \def\GTU {\mathfrak U}    \def\GTV {\mathfrak V}
\def\GTW {\mathfrak W}   \def\GTX {\mathfrak X}    \def\GTY {\mathfrak Y}
\def\GTZ {\mathfrak Z}   \def\GTQ {\mathfrak Q}

\font\Sym= msam10 
\def\SYM#1{\hbox{\Sym #1}}
\newcommand{\bdw}{\prt\Gw\xspace}

\date{}
\maketitle\medskip

\noindent{\small {\bf Abstract} We study the semilinear elliptic equation $-\Gd u + g(u)\gs = \gm$ with Dirichlet boundary condition in a smooth bounded domain where $\gs$ is a nonnegative Radon measure,  $\gm$ a Radon measure and $g$ is an absorbing nonlinearity. 
We show that the problem is well posed if we assume that $\gs$ belongs to some Morrey class. Under this condition we give a general existence result for any bounded measure provided $g$ satisfies a subcritical integral assumption. We study also the supercritical case when 
$g(r)= \abs r^{q-1}r$, with $q>1$ and $\gm$ satisfies an absolute continuity condition expressed in terms of some capacities involving $\gs$. 
\medskip

\vspace{1mm}
\hspace{.05in}
\tableofcontents
\medskip
\mysection{Introduction}
Let  $\Gw \sbs \BBR^N$ be a bounded domain with a $C^2$ boundary, $\gs$ a nonnegative Radon measure in $\Gw$ and $g:\BBR \to\BBR$ a continuous function satisfying, for some $r_0\geq 0$,
\bel{Z0}
rg(r)\geq 0 \qquad\text{for all } r\in (-\infty,-r_0]\cup[r_0,\infty).
\ee 
In this article we consider the following problem
\bel{Z1}\BA {lll} -\Gd u + g(u)\gs = \gm \qq &\text{in } \Gw\\
\phantom{-\Gd  + g(u)\gs}
u=0\qq &\text{in } \prt\Gw,
\EA\ee
where $\gm$ is a Radon measure defined in $\Gw$. By a solution we mean a function $u\in L^1(\Gw)$ such that $\gr g(u)\in L_\gs^1(\Gw)$, where  $\gr(x)=\dist(x,\prt\Gw)$
and $L_\gs^1(\Gw)$ is the Lebesgue space of functions integrable with respect to $\sigma$, satisfying 
\bel{Z2}\BA {lll} 
-\myint{\Gw}{}u\Gd\gz dx+\myint{\Gw}{}g(u)\gz d\gs=\myint{\Gw}{}\gz d\gm,
\EA\ee
for all $\gz\in W^{1,\infty}_0(\Gw)$ such that $\Gd\gz\in L^\infty(\Gw)$. In the sequel, such a solution is called a {\it very weak solution}. A measure $\gm$ such that the problem admits a solution is called a {\it good measure}. We emphasize on the particular cases where $g(r)=\abs r^{q-1}r$ with $q>0$, or $g(r)=e^{\ga r}-1$ with $\ga>0$ and $N= 2$.

  \smallskip

When $\gs$ is a measure with constant positive density with respect to the Lebesgue measure in $\BBR^N$, this problem has been initiated by Brezis and Benilan \cite{BrBe0}, \cite{BrBe1} who gave a general existence result for any bounded measure $\gm$ under an integrability condition of $g$ at infinity; their proof is based upon an a priori estimate of approximate solutions $u_n$ in Lorentz spaces $L^{q,\infty}(\Gw)$, yielding the uniform integrability of $g(u_n)$ and hence the pre-compactness in $L^1(\Gw)$. If $g(r)=\abs r^{q-1}r$, integrability condition is fufilled if and only if $0<q<\frac{N}{N-2}$ (any $q>0$ if $N=2$). In the 2-dim case the integrability condition have been replaced by the exponential order of growth of $g$ in \cite{Va1}. When 
 $g(u)=\abs{u}^{q-1}u$ with $q\geq \frac{N}{N-2}$ not any bounded measure is eligible for solving $(\ref{Z1})$. In fact Baras and Pierre \cite{BaPi1} proved that when $N>2$ and $q>1$, a bounded Radon measure $\gm$ is eligible if and only if it vanishes on Borel sets with $c_{2,q'}$-capacity zero, where  $q'=\frac{q}{q-1}$ is the conjugate exponent of $q$. 
Contrary to the previous subcritical case, the method for proving the necessity of this condition is based upon a duality-convexity argument, while the sufficiency uses the fact that any  positive Radon measure absolutely continuous with respect to the  $c_{2,q'}$-capacity can be approximated from below by an non-decreasing sequence of positive measures in 
$W^{-2,q}(\Gw)$ (see \cite{FePr}). 
Furthermore they also gave  a necessary and sufficient condition for a compact subset $K\subset\Gw$ to be removable for equation 
  \bel{Z3}\BA {lll} 
-\Gd u+\abs u^{q-1}u=0\qquad\text{in }\,\Gw\setminus K,
\EA\ee
 namely that $c_{2,q'}(K)=0$.  \smallskip
 
 The aim of this paper is to extend the previous constructions of Benilan-Brezis, Baras-Pierre and Vazquez to the case where $\gs$ is a general measure.  In order to be able to deal with the convergence of approximate solutions we assume that $\gs$ belongs to the Morrey class $\CM^+^+_{\frac  N{N-\gth}}(\Gw)$ for some  $\theta\in [0,N]$} which means
   \bel{Z4}\BA {lll} 
\abs{B_r(x)}_\gs:=\myint{B_r(x)}{}d\gs\leq cr^{\gth}\qquad\text{for all } (x,r)\in \Gw\ti (0,\infty),
\EA\ee
for some $c>0$. Note that we extend $\gs$ by  0  in $\BBR^N\setminus\Gw$. 

  \smallskip

Our first result is the following:\medskip

\nind{\bf Theorem A}  {\it Assume $\gs\in \CM^+^+_{\frac  N{N-\gth}}(\Gw)$ for some  
$ \theta\in (N-2,N]$ and that $g$ satisfies $(\ref{Z0})$. Then, for any $\gm\in L_\gr^1(\Gw)$, there exists a very weak solution $u$ of problem $(\ref{Z2})$. 
If we assume moreover that $g$ is nondecreasing and if $u'$ is a very weak solution of $(\ref{Z2})$ with right-hand side 
$\gm'\in L_\gr^1(\Gw)$, then the following estimates hold
   \bel{Z5}\BA {lll} 
-\myint{\Gw}{}\abs{u-u'}\Gd\gz dx+\myint{\Gw}{}\abs{g(u)-g(u')}\gz d\gs
\leq \myint{\Gw}{}\abs{\gm-\gm'}dx, 
\EA\ee
and
   \bel{Z6}\BA {lll} 
-\myint{\Gw}{}(u-u')_+\Gd\gz dx+\myint{\Gw}{}(g(u)-g(u'))_+\gz d\gs\leq \myint{\Gw}{}(\gm-\gm')_+dx
\EA\ee
for all $\gz\in W^{1,\infty}_0(\Gw)$  such that  $\Gd\gz\in L^\infty(\Gw)$ and 
$\gz\geq 0$.}

\medskip  

\noindent Note that $(\ref{Z5})$ implies the uniqueness of the solution of $(\ref{Z2})$, that we denote by $u_\gm$, and $(\ref{Z6})$ the monotonicity of the mapping $\gm\mapsto u_\gm$. 

\medskip

The next result extends Benilan-Brezis unconditional existence result for measures.\medskip

\nind{\bf Theorem B}  {\it Let $N>2$ and $\gs\in \CM^+^+_{\frac  N{N-\gth}}(\Gw)$ with $N\geq \gth> N-\frac{N}{N-1}$. Assume that $g$ satisfies $(\ref{Z0})$ and  
$\abs {g(r)}\leq \tilde g(\abs r)$ for all $\abs r\geq r_0$ where $\tilde g$ is a continuous nondecreasing function on $[r_0,\infty)$ verifying 
   \bel{Z7}\BA {lll} 
\myint{r_0}{\infty}\tilde g(t)t^{-1-\frac{\gth}{N-2}}dt<\infty.
\EA\ee
Then, for any bounded Radon measure $\gm$, there exists a very weak solution $u$ of problem $(\ref{Z2})$ which moreover belongs to $L^1_\sigma(\Omega)$. 
Moreover, if we assume that $g$ is nondecreasing then the solution is unique.}

\medskip 

\noindent Note that we recover Benilan-Brezis result when $\sigma$ is the Lebesgue mesure (so that $\theta=N$). 
Note also that when $g(r)=\abs r^{q-1}r$, the integrability condition $(\ref{Z7})$ is fullfilled if and only if $0<q<\frac{\theta}{N-2}$.

\medskip

In the 2-dimensional case the condition on $\gth$ is $2\geq \gth>0$ but $(\ref{Z7})$ has to be modified. If $ f:\BBR\mapsto\BBR_+$ is nondecreasing we define its exponential order of growth at $\infty$ (see \cite{Va1}) by
   \bel{Z9}\BA {lll} 
a_\infty(f)=\inf\left\{\ga\geq 0:\myint{0}{\infty}f(s)e^{-\ga s}ds<\infty\right\}.
\EA\ee
Similarly, if $h:\BBR\mapsto\BBR_-$ is nondecreasing its exponential order of growth at $-\infty$ is 
   \bel{Z9'}\BA {lll} 
a_{\text{-}\infty}(h)=\sup\left\{\ga\leq 0:\myint{-\infty}{0}h(s)e^{\ga s}ds >-\infty\right\}.
   \EA\ee

If $g:\BBR\mapsto\BBR$ satisfies $(\ref{Z0})$ but is not necessarily nondecreasing, we define the monotone nondecreasing hull $g^*$  of $g$ by
   \bel{Z'9}\BA {lll} 
g^*(r)=\left\{\BA {lll}\sup\{g(s):s\leq r\}\qquad&\text{for all } r\geq r_0\\
0\qquad&\text{for all } r\in (-r_0, r_0)\\
\inf\{g(s):s\geq r\}\qquad&\text{for all } r\leq -r_0.\EA\right.
\EA\ee
We set
   \bel{Z''9}\BA {lll} 
  a_\infty(g)= a_\infty(g^*_+)\quad\text{and }\;a_{\text{-}\infty}(g)=a_{\text{-}\infty}(g^*_-).
\EA\ee\smallskip

\nind{\bf Theorem C}  {\it Let $\gs\in \CM^+^+_{\frac  2{2-\gth}}(\Gw)$ with $2\geq \gth>0$ and $g:\BBR\mapsto\BBR$ satisfies $(\ref{Z0})$. \smallskip

\nind (I) If $  a_\infty(g)=0=a_{\text{-}\infty}(g)$, then for any $\gm\in\mathfrak M_b(\Gw)$, problem $(\ref{Z2})$ admits a very weak solution.\smallskip

\nind (II) If $0<  a_\infty(g)<\infty$ and $-\infty<a_{\text{-}\infty}(g)<0$ there exists $\gd>0$ such that if $\gm\in\mathfrak M_b(\Gw)$ satisfies 
$\norm\gm_{\mathfrak M_b}\leq\gd$ problem $(\ref{Z2})$ admits a very weak solution.
}

\medskip

In the {\it supercritical case}, that is when $(\ref{Z7})$ is not satisfied, all the measures are not eligible for solving $(\ref{Z2})$. Following \cite{MaVe1}, \cite[Th 4.2 ]{Ve0} we can give a sufficient existence condition involving the Green function of the Laplacian. 
Let $G(.,.)$ be the Green kernel defined in $\Gw\ti\Gw$ and $\BBG[.]$ the corresponding potential operator acting on bounded measures $\nu$ namely 
$\BBG[\nu](x)=\int_\Omega G(x,y)\,d\nu(y)$. We have the following result: 

\medskip

\nind{\bf Theorem D}  {\it Let $\gs\in \CM^+^+_{\frac  N{N-\gth}}(\Gw)$ with $N\geq \gth>N-\frac{N}{N-1}$ and assume that $g$ is nondecreasing and vanishes at $0$. \smallskip

\nind (I) If 
$\gm\in \mathfrak M_b(\Gw)$ satisfies
\bel{Z11}\BA {lll} 
\gr g(\BBG[\abs{\gm}])\in L^1_\gs(\Gw),
\EA\ee
then problem $(\ref{Z2})$ admits a unique very weak solution.\smallskip

\nind (II) Let $\gm=\gm_r+\gm_s$ where $\gm_r$ is absolutely continuous with respect to the Lebesgue measure and $\gm_s$ is singular. Assume that $g$ satisfies the $\Gd_2$ condition, namely that 
\bel{Z12}\BA {lll} 
\abs{g(r+r')}\leq a\left(\abs{g(r)}+\abs{g(r')}\right)+b\qquad\text{for all } r,r'\in\BBR,
\EA\ee
for some $a>1$ and $b\geq 0$. Then the previous assertion holds if $(\ref{Z11})$ is replaced by
\bel{Z13}\BA {lll} 
\gr g(\BBG[\abs{\gm_s}])\in L^1_\gs(\Gw).
\EA\ee
}

\medskip

\noindent 
Notice that $(\ref{Z11})$  holds if either (i) $\sigma$ and $\lambda$ have disjoint support, or (ii) $\gm\in\CM^+_p(\Gw)$ for some $p>\frac{N}{2}$.  
Indeed if (i) holds then $\BBG[\abs{\gm}]$ is bounded pointwise on the support of $\sigma$, and if (ii) holds then by \rlemma{Miy} $\BBG[\abs{\gm}]$ is bounded pointwise in $\Omega$. Obviously the same comment holds in the setting of II.

\medskip

In order to make more explicit conditions $(\ref{Z11})$, $(\ref{Z13})$, we introduce the following growth assumption on $g$: 
\bel{Z14}\BA {lll} 
\abs{g(r)}\leq c(1+\abs r^q)\qquad\text{for all } r\in\BBR,
\EA\ee
for some $q>1$. 
Notice that $\tilde g(r)=1+r^q$ satisfies $(\ref{Z7})$ if and only if $q < \frac{\gth}{N-2}$. 
When $\sigma$ is the Lebesgue measure and $g(r)=|r|^{q-1}r$, Baras and Pierre \cite{BaPi1} gave a necessary and sufficient condition for the existence of a solution to $(\ref{Z1})$ involving certain capacity associated to the Bessel potential spaces $H^{s,p}(\BBR^N)$ where $s\in\BBR$ and $p\in [1,\infty]$. Let us recall that 
   \bel{Z15}\BA {lll} 
H^{s,p}(\BBR^N)=\left\{f: f={\bf G}_s\ast h, h\in L^p(\BBR^N)\right\},
\EA\ee
where ${\bf G}_s$ is the Bessel kernel of order $s$. By extension ${\bf G}_0=\gd_0$, hence $H^{s,p}(\BBR^N)=L^p(\BBR^N)$.
When $s$ is a positive integer, it is proved by Calder\'on \cite[Theorem 1.2.3]{AdHe} that $H^{s,p}(\BBR^N)$ is the standard Sobolev space 
$W^{s,p}(\BBR^N)$.  
If $s>0$, we denote by $c_{s,p}$ the associated capacity, called the Bessel capacity. It is defined for any compact set $K\subset \mathbb{R}^N$ by 
\begin{equation}\label{DefCapBessel} 
  c_{s,p}(K) = \inf\,\{ \|\phi\|_{H^{s,p}}^p\,:\, \phi\in\mathcal{S}(\mathbb{R}^N),\,
\phi\geq 1 \mbox{ on } K\}. 
\end{equation} 
The definition of $c_{s,p}$ is then extended first to open sets and then to arbitrary sets. 
We refer to \cite{AdHe} for general properties of Bessel spaces and their associated capacities 
$c_{s,p}$. 
We say that a measure  $\gm\in\mathfrak M_b(\Gw)$ is absolutely continuous  with respect to the $c_{s,p}$-capacity if for any Borel subset $E\subset \mathbb{R}^N$, 
$$ c_{s,p}(E)=0 \Longrightarrow |\mu|(E)=0. $$ 
Baras and Pierre's result states that equation $(\ref{Z1})$, with $\sigma$ standing for the Lebesgue measure and $g(r)=|r|^{q-1}r$, has a solution if and only if $\mu$ is absolutely continuous with respect to the $c_{2,q'}$-capacity. The next result generalizes the "if" part to the case where $\sigma$ belongs to some Morrey space.

\medskip

\nind{\bf Theorem E}  {\it Let $\gs\in \CM^+^+_{\frac  N{N-\gth}}(\Gw)$ with $N\geq \gth> N-\frac{N}{N-1}$ and assume that $g$ is nondecreasing and satisfies $(\ref{Z0})$ and $(\ref{Z14})$. Let  $p>1$ and $s\geq 0$  such that $N>sp>N-\theta$ and $\frac{\gth p}{N-sp}\geq q$.  If $\gm\in\mathfrak M_b(\Gw)$ is absolutely continuous  with respect to the $c_{2-s,p'}$-capacity, then $(\ref{Z1})$ admits a unique very weak solution.}

\medskip

\noindent As a particular case, we take $p=q$ and obtain that if $\gm$ is absolutely continuous  with respect to the $c_{2-\frac{N-\gth}{q},q'}$-capacity, then $(\ref{Z2})$ admits a unique  solution. 
We thus recover Baras-Pierre's sufficient condition \cite{BaPi1} when $\theta=N$. 

\medskip 

\noindent 
We give an explicit condition on the measure $\gm$ in terms of Morrey spaces implying that it satisfies the conditions of Theorem E.

\bprop{CondSimple} Under the assumptions on $\sigma$ and $g$ of Theorem E, if $\mu\in \CM^+_{\frac  N{N-\gth^*}}(\Gw)$
for some $\gth^*>\frac{(N-2)q-\theta}{q-1}$, then $(\ref{Z2})$ admits a unique very weak solution.
\es

\smallskip 

\noindent Notice that the condition on $\mu$ given in \rprop{CondSimple} is weaker than the one given after Theorem D.

\medskip

When $g(r)=\abs r^{q-1}r$ with $q>1$, one can find a necessary conditions for the existence 
 of a solution of $(\ref{Z2})$ in the supercritical case under additional regularity assumptions on $\gs$. By \cite[Def 2.3.3, Prop. 2.3.5]{AdHe}, the following expression 
   \bel{Z16}\BA {lll} 
c_{q}^\gs(E)=\inf\left\{\myint{\Gw}{}\abs v^{q'}d\gs:v\in L^{q'}_\gs(\Gw),\,v\geq 0, \,\BBG[v\gs]\geq 1 \mbox{ on $E$}\right\},
\EA\ee
where $E$ is any subset of $\Gw$ defines an outer capacity. 
The measure is called $\gth$-regular if
$$\myfrac{1}{c}r^\gth\leq\myint{B_r(x)}{}d\gs\leq c r^\gth\qquad \text{for all } (x,r)\in \Gw\ti (0,1],
$$
The next result gives a necessary condition for a measure to be a good measure. \medskip 

\nind{\bf Theorem F}  {\it Let $q>1$ and $\gs\in \CM^+^+_{\frac  N{N-\gth}}(\Gw)$ be $\gth$-regular with $N\geq \gth>N-2$. 
If $\gm\in\mathfrak M^+_b(\Gw) $ is such that problem $(\ref{Z2})$ with $g(r)=\abs r^{q-1}r$ admits a very weak solution, then 
$\gm$ vanishes on any Borel set $E$ such that  $c_{q}^\gs(E)=0$.}

\medskip

Furthermore the $c_{q}^\gs$- capacity admits the following representation in terms of Besov capacities. 
If  $\Gg\subset\Gw$ is the support of $\gs$,  we denote by $B^{2-\frac{N-\gth}{q},\Gg}_{q',\infty}(\Gw)$ the closed subspace of  distributions $\gz\in B^{2-\frac{N-\gth}{q}}_{q',\infty}(\Gw)$  such that the support of the distribution $\Gd\gz$ is a subset of $\Gg$. Then 
\bel{Z18}\BA {lll} 
c_{q}^\gs(K)\sim c^{2-\frac{N-\gth}{q},\Gg}_{q',\infty}(K):=\inf\left\{\norm\gz^{q'}_{B^{2-\frac{N-\gth}{q}}_{q',\infty}}: \gz\in B^{2-\frac{N-\gth}{q},\Gg}_{q',\infty}(\Gw),\,\gz\geq \chi_{_K}\right\},
\EA\ee
for all compact subset $K\subset\Gw$.\medskip

Finally a complete characterization of removable sets can be obtained under a much stronger assumption on $\gs$, namely that $d\gs=w dx$ with $\gw:=w^{-\frac{1}{q-1}}\in L^1_{loc}(\Gw)$. If $K\subset\Gw$ is compact, we set

 \bel{Z19}\BA {lll} 
c^\gw_{q}(K)=\inf\left\{\myint{\Gw}{}\abs{\Gd \gz}^{q'}\gw dx: \gz\in C^{\infty}_0(\Gw), 0\leq \gz\leq 1,\gz=1\text{ in a neighborhood of }K\right\}.
\EA\ee
This defines a capacity on Borel sets of $\Gw$. \medskip

\nind{\bf Theorem G}. {\it Assume $q>1$ and there exists a nonnegative Borel function $w$ in $\Gw$ in the Muckenhoupt class  $A_q(\Gw)$ such that $d\gs=w dx$. If $K\subset\Gw$ is compact, a function $u\in L^1_{loc}(\Gw\setminus K)$ such that $|u|^qw\in L^1_{loc}(\Gw\setminus K)$ which satisfies 
 \bel{Z20}\BA {lll} 
-\Gd u+w\abs u^{q-1}u=0,
\EA\ee
in the sense  of distributions in $\Gw\setminus K$ can be extended as a solution of the same equation in whole $\Gw$ if and only if 
$c_{q,w}(K)=0$.
 }
 
The assumption $w\in A_q(\Gw)$ can be weakened and replaced by $\gw=w^{\frac{1}{1-q}}$ is $q'$-admissible in the sense of \cite[Chap 1]{HKM}, a condition which implies in particular the validity of the Gagliardo-Nirenberg and the Poincar\'e inequalities. 
 
 \medskip

\mysection{Preliminaries}

In the whole paper $c$ denotes a generic positive constant whose  value can change from one ocurrence to another even within a single string of estimates.
Sometimes, in order to avoid ambiguity, we are led to introduce other notations for  constant, for example  $c'$. 

\medskip 

We denote by $\mathfrak M_b(\Gw)$ the space of outer regular bounded Borel measures on $\Gw$ equipped with the total variation norm, and by $\mathfrak M_b^+(\Gw)$ its positive cone. 
Since $\Gw$ is bounded we can identify bounded Radon measures in $\Gw$ with measures $\gm$ in $\overline\Gw$ such that $\abs\gm(\prt\Gw)=0$.
All the measures are extended by 0 in $\mathbb{R}^N\backslash \Omega$. 

\medskip

Let $G(.,.)$ be the Green kernel defined in $\Gw\ti\Gw$ and $\BBG[.]$ the corresponding potential operator acting on bounded measures $\nu$ namely 
$\BBG[\nu](x)=\int_\Omega G(x,y)\,d\nu(y)$.
We denote  $L^{p,\infty}(\Gw)$ the  usual weak $L^p$ space. 
The next result is classical and valid in a much more general setting (see e.g. \cite{BBC}, \cite {DHM}).

\blemma{L1} Let $\gm\in\mathfrak M_b(\Gw)$ and $v=\BBG[\gm]$ be the (very weak) solution of
 \bel{Y2+2}\BA {lll} \displaystyle
-\Gd v=\gm\qquad&\text{in }\;\Gw\\
\phantom{-\Gd}
v=0&\text{in }\;\prt\Gw.
\EA\ee
I- If $N\geq 2$, then $v\in L^{\frac {N}{N-2},\infty}(\Gw)$, $\nabla v\in L^{\frac {N}{N-1},\infty}(\Gw)$ and 
\bel{Y2-1}
\norm v_{L^{\frac {N}{N-2},\infty}}+\norm {\nabla v}_{L^{\frac {N}{N-1},\infty}}\leq c\norm\gm_{\mathfrak M_b}.
\ee
II- If $N= 2$, then $v\in BMO(\Gw)$, $\nabla v\in L^{2,\infty}(\Gw)$ and
\bel{Y2-2}
\norm v_{BMO}+\norm {\nabla v}_{L^{2,\infty}}\leq c\norm\gm_{\mathfrak M_b}.
\ee
\es

This result can be refined when more information is available on the degree of concentration of $\mu$. This lead to the definition of Morrey spaces of measures.

\subsection{Morrey spaces of measures}

If $1\leq p\leq\infty$ we define the Morrey space 
$\CM^+ _{p}(\Gw)$ as the set of bounded outer regular Borel measures $\gm$ defined in $\Gw$ and extended by $0$ in $\Gw^c$, satisfying
\bel{Y1}
\abs {B_r(x)}_{\gm}:=\myint{B_r(x)}{}d\abs\gm\leq cr^{N(1-\frac1p)}\qquad\text{for all } (x,r)\in \Gw\ti\BBR_+,
\ee
for some $c>0$. In particular $\mu\in \CM^+ _\frac{N}{N-\theta}(\Gw)$, $\theta\in [0,N]$, if 
$$\myint{B_r(x)}{}d\abs\gm\leq cr^\theta \qquad\text{for all } (x,r)\in \Gw\ti\BBR_+.$$ 
We refer to \cite{Mi} for a detailed study of $\CM^+ _{p}(\Gw)$ and full proofs of the various results we will recall now. Endowed with the norm
\bel{Y2}
\norm\gm_{\CM^+_{p}}=\sup_{ (x,r)\in \Gw\ti\BBR_+}r^{N(\frac1p-1)}\abs {B_r(x)}_{\gm},
\ee
$\CM^+_{p}(\Gw)$ is a Banach space and $\CM^+^+_p(\Gw)=\CM^+_{p}(\Gw)\cap\mathfrak M_b^+(\Gw)$ is its positive cone. We also set $M_{p}(\Gw)=\CM^+_{p}(\Gw)\cap L^1_{loc}(\Gw)$; it is a closed subspace of $\CM^+_{p}(\Gw)$ and,  if $1<p<\infty$,  the following imbedding holds
\bel{Y2+1}
L^{p}(\Gw)\hookrightarrow L^{p,\infty}(\Gw)\hookrightarrow M_{p}(\Gw). 
\ee

Note that since $\Omega$ is bounded and any measure in $\Omega$ is extended to $\mathbb{R}^N$ by 0,  it is easily seen that if $1\leq q\leq p\leq \infty$ we have a continuous embedding $\CM^+_p(\Gw)\hookrightarrow \CM^+_q(\Omega)$ with 
\bel{Y7}
\norm {v}_{\CM^+_{q}}\leq ({\rm diam}(\Gw))^{\frac{N}{q}-\frac{N}{p}}\norm {v}_{\CM^+_{p}}\qquad\text{for all } v\in \CM^+_p(\Gw).
\ee
Indeed for any $x\in\Omega$ the ball centered at $x$ with radius ${\rm diam}(\Gw)$ contains $\Omega$ so that it is enough to consider $r\leq {\rm diam}(\Gw)$. We have 
$$ r^{-N(1-1/q)}\abs{B_r(x)}_{\gm}
\leq r^{-N(1-1/q)} \norm\gm_{\CM^+_p} r^{N(1-1/p)}
\leq ({\rm diam}(\Gw))^{\frac{N}{q}-\frac{N}{p}}     \norm\gm_{\CM^+_p}. $$ 

The following imbedding inequalities holds.

\blemma{Miy} Let $\gm\in\CM^+_{p}(\Gw)$ and $v$ be the solution of $(\ref{Y2+2})$.\smallskip

\nind I- If $1<p<\frac N2$, then $v\in M_q(\Omega)$ 
with $\frac{1}{q}=\frac{1}{p}-\frac{2}{N}$ and there holds 
\bel{Y2-3}
\norm v_ {\CM^+_{q}}    \leq c\norm\gm_{\CM^+_p}.
\ee
\nind II- If $p>\frac N2$, then $v$ is bounded pointwise 
and
\bel{Y2+3}\BA {lll} \displaystyle
(i)\;\;\qquad\phantom{,\cap C^{\ga}(\BBR^N)}   
v(x) \leq c\norm\gm_{\CM^+_p}\qquad\text{for all }x\in \Gw,\\[2mm]
(ii)\phantom{i}\displaystyle\;\;\qquad\sup_{x\neq y}\myfrac{\abs{v(x)-v(y)}}{\abs{x-y}^\ga}\leq c\norm\gm_{\CM^+_{p}}\,\text{ with }\;\ga=2-\frac{N}{p}\,&\text{ if }\;N>p>\myfrac{N}{2},\\[3mm]
(iii)\displaystyle\;\;\qquad\sup_{x\neq y}\myfrac{\abs{v(x)-v(y)}}{\abs{x-y}^\ga}\leq c\norm\gm_{\CM^+_{p}}\,\text{ with }\;\ga\in (0,1)\,&\text{ if }\;N=p,\\[3mm]
(iv)\displaystyle\;\;\phantom{---}\qquad\sup_{x}\abs{\nabla v(x)}\leq c\norm\gm_{\CM^+_{p}}&\text{ if }\;N<p.
\EA\ee
\es

\nind\Remark 
The previous regularity results are proved in \cite[Prop. 3.1, 3.5]{Mi} when $v=I_\alpha*\mu$ where $I_\alpha$ is the Riesz potential. However it is easily seen that the proof in  \cite{Mi} can be adapted to our setting. 
In particular for $(\ref{Y2-3})$ we need that $G(x,y)\le c|x-y|^{2-N}$, for (i) we use  $(\ref{Y7})$. 
\medskip

\nind\Remark If we assume that $\gm\in\mathfrak M_\gr(\Gw)\cap \CM^+_{p,loc}(\Gw)$, the previous estimates acquire a local aspect and remain valid 
provided the supremum in the norms on the left-hand sides are taken  on compact subsets of $\Gw$.\smallskip

\subsection{Trace embeddings}

Some applications of Morrey spaces to imbedding theorems (also called trace inequalities) can be found in Adams-Hedberg's book \cite{AdHe}. For the sake of completeness, we quote here the main result therein we will use in the sequel. If $0<\ga<N$ we recall that $I_{\ga}$ (resp. $G_\ga$) is the Riesz potential (resp. the Bessel potential) of order $\ga$ in $\BBR^N$. The next result is \cite[Th 7.2.2, 7.3.2 ]{AdHe} (recall that the $c_{_{I_\ga},p}\,$-Riesz capacity of a ball $B_r(x)$ is proportional to $r^{N-\alpha p}$ - see \cite[Prop. 5.1.2]{AdHe}.)

\bprop{AH} Let  $\gs$ be a nonnegative Radon measure in $\BBR^N$, $N>\ga p$ and $1<p<q<\frac{Np}{N-\ga p}$.\smallskip

\nind (I)-  The following assertions are equivalent:
\bel{Y3}
\norm{I_\ga\ast f}_{L_\gs^{q}(\BBR^N)}\leq c_1\norm {f}_{L^p(\BBR^N)}\qquad\text{for all } f\in L^p(\BBR^N),
\ee
for some $c_1=c_1(N,\ga,p,q)>0$, and
\bel{Y3'}
\gs\in \CM^+_{r}(\BBR^N)\quad\text{with }\, \frac{1}{r}=q\left(\frac{1}{q}-\frac{1}{p}+\frac{\ga}{N}\right).
\ee

\nind (II)- The mapping $f\mapsto G_\ga\ast f$ is continuous from $L^p(\BBR^N)$ to $L_\gs^{q}(\BBR^N)$ if and only if 
\bel{Y3'''}
\gs(K)^\frac{1}{q}\leq c_2 \left(c_{\ga,p}(K)\right)^\frac{1}{p}  \qquad \text{for all }\, K\subset \mathbb{R}^N,
\ee
where $c_{\ga,p}$ denotes the Bessel capacity of order $\ga$ defined in $(\ref{DefCapBessel})$. In fact this holds if and only if 
\bel{Y3''}
\gs(B_r(x))\leq c_3 \left(c_{\ga,p}(B_r(x))\right)^{q/p}  \qquad \text{for all }\, x\in\mathbb{R}^N,\,0<r\leq 1.  
\ee

\nind (III)-  A necessary and sufficient condition in order the mapping $f\mapsto G_\ga\ast f$ be compact from $L^p(\BBR^N)$ to $L_\gs^{q}(\BBR^N)$ is
\bel{Y4}\BA{lll}\displaystyle
(i)\qquad\lim_{\gd\to 0}\sup_{x\in\BBR^N\!,\,r
\leq\gd}\myfrac{\gs(B_r(x))}{\left(c_{\ga,p}(B_r(x))\right)^{\frac qp}}=0
\\[6mm]\displaystyle
(ii)\qquad\lim_{\abs x\to \infty}\sup_{r\leq 1}\myfrac{\gs(B_r(x))}{\left(c_{\ga,p}(B_r(x))\right)^{\frac qp}}=0. 
\EA\ee
\es

If $\BBR^N$ is replaced by a smooth bounded set $\Gw$, we extend  any bounded Radon measure in $\Gw$ by zero in $\Gw^c$. 
In view of \cite[5.6.1] {AdHe} the $c_{_{I_\ga},p}\,$-Riesz capacity and $c_{\ga,p}\,$- Bessel capacity of balls $B_r(x)$ with $x\in\Omega$ and $r\le 1$ are then equivalent. It follows that 
$c_{\ga,p}(B_r(x))\simeq r^{N-\alpha p}$. 
Then, it follows from II and III above, the definition of $H^{\alpha,p}(\mathbb{R}^N)$ and the existence of an extension operator 
$H^{\alpha,p}(\Omega)\hookrightarrow H^{\alpha,p}(\mathbb{R}^N)$  that the following holds,

\bprop{AdamsMazya} Under the assumptions of \rprop{AH}, 
the embedding $H^{\alpha,p}(\Omega)\hookrightarrow L_\gs^{q}(\Omega)$ is: 

\nind (I)-  continuous if and only if 
$\left(\gs(K)\right)^\frac{1}{q}\leq c_2 \left(c_{\ga,p}(K)\right)^\frac{1}{p}$ for all $K\subset \mathbb{R}^N$, i.e. if and only if  $\gs\in \CM^+^+_{r}(\BBR^N)$ with 
$\frac{1}{r}=q\left(\frac{1}{q}-\frac{1}{p}+\frac{\ga}{N}\right)$.

\nind (II)-  compact if and only if 
\bel{Y8}
\lim_{r\to 0}\sup_{x\in\Gw}\myfrac{\gs(B_r(x))}{r^{\frac{(N-\ga p)q}{p}}}=0.
\ee
\es

\noindent As an immediate corollary, 

\bprop{AdamsMazya2} 
Let  $\gs\in \CM^+^+_\frac{N}{N-\theta}(\Omega)$, i.e. $\sigma(B_r(x))\leq cr^\theta$, $N>\ga p$ and $1<p<q<\frac{Np}{N-\ga p}$.  Then the embedding 
\bel{Y31}
H^{\alpha,p}(\Omega)\hookrightarrow L_\gs^{q}(\Omega), 
\ee
is continuous iff $\gs(K)\leq c_1 \left(c_{\ga,p}(K)\right)^\frac{q}{p}$ for all $K\subset \mathbb{R}^N$ which holds iff 
 $q\leq \frac{\theta p}{N-\alpha p}$. And the embedding $(\ref{Y31})$ is compact iff 
$q< \frac{\theta p}{N-\alpha p}$.
\es

\medskip

Other trace inequalities can be found in \cite{MazVer}. In the case $N=\ga p$ the following estimate holds, see e.g. \cite{Ad1}, \cite[Corollary 8.6.2]{Maz1}, \cite{Yu}.
\bprop{AM} Let  $\gs$ be a nonnegative Radon measure in $\BBR^N$ with compact support and $N=\ga p$, $p>1$. Then there exists a constant $b=b(N,\ga,p)>0$ such that 
\bel{Y5}\displaystyle
\sup_{\norm f_{L^p}\leq 1}\myint{\BBR^N}{}\exp\left(b\abs{G_\ga\ast f}^{p'}\right)d\gs<\infty
\ee
if and only if $\gs\in\CM^+^+_\gt(\BBR^N)$ for some $\gt\in (1,\infty)$.
\es

When $p=1$ the next result is proved in \cite[Sec 1.4.3]{Maz1}

\bprop{MaSh} Let  $\gs$ be a nonnegative bounded Radon measure in $\BBR^N$, $\ga$ be an integer such that $1\leq \ga\leq N$ and $q\geq 1$. Then the following estimate holds  
\bel{Y6}\displaystyle
\norm f_{L^q_\gs}\leq c_2\sum_{\abs\beta=\ga}\|D^\alpha f\|_1 
\qquad\text{for all } f\in C^\infty_0(\BBR^N),
\ee
for some $c_2=c_2(N,p,q,\ga)>0$ if and only if 
$\gs\in \CM^+^+_{\frac{N}{N - q(N-\alpha)}}(\BBR^N)$.
\es

\mysection{The subcritical case}

\subsection{The variational construction}

We prove in this section that if $\gm\in W^{-1,2}(\Gw)$ then, under some assumptions on $g$ and $\sigma$, equation $(\ref{Z1})$ has a variational solution. 

We assume that $g\in C(\BBR)$ satisfies $(\ref{Z0})$, and set $G(r):=\myint{0}{r}g(s)ds$. 
We will find a solution to $(\ref{Z1})$ minimizing the functional 
\bel{Y10}
J(v):=\myfrac{1}{2}\myint{\Gw}{}\abs{\nabla v}^2dx+\myint{\Gw}{}G(v)\,d\gs-\langle\gm, v\rangle,
\ee
over the set 
\bel{Y9}
X_G(\Gw):=\{v\in W^{1,2}_0(\Gw):G(v)\in L^1_\gs(\Gw)\}.
\ee
The next proposition is a variant of a result in \cite{BrBr}.
\bprop{varia} Assume $\gs\in\CM^+^+_{\frac N{N-\gth}}(\Gw)$ with $N\geq \gth>\frac N2-1$. 
If $\gm\in W^{-1,2}(\Gw)$ there exists $u\in X_G(\Gw)$ which minimizes $J$ in $X_G(\Gw)$. Furthermore $u$ is a weak solution of $(\ref{Z1})$ in the sense that
\bel{Y10-1}
\myint{\Gw}{}\nabla u.\nabla\gz dx+\myint{\Gw}{}g(u)\gz d\gs=\langle\gm,\gz\rangle\quad\text{for all }\;\gz\in C^\infty_0(\Gw).
\ee
If $g$ is nondecreasing this solution is unique and denoted by $u_\gm$, and the mapping $\gm\mapsto u_\gm$ is nonnecreasing. 
\es
\noindent\Proof {\it Step 1: Existence of a minimizer}. If $N>2$ we apply $(\ref{Y31})$  with $\ga=1$ and $p=2$, recalling that by Fourier transform  $H^{1,2}(\Omega)=W^{1,2}(\Omega)$ (it is a special case of Calder\'on's theorem), to obtain that 
\bel{Y11} W^{1,2}_0(\Gw) \hookrightarrow      L^{\frac{2\gth}{N-2}}_\gs(\Gw). 
\ee
If $N=2$ with $p=2$ we take any $\ga<1$ and obtain
\bel{Y11'} 
\norm{f}_{L^{\frac{\gth}{1-\ga}}_\gs}\leq c_1\norm{f}_{W^{\ga,2}}\leq c'_1\norm{f}_{W^{1,2}}. 
\ee
According to  \rprop{AdamsMazya2} the imbedding of $W^{1,2}_0(\Gw)$ into $L^p_\gs(\Gw)$ is compact for any $p\in [1,\frac{2\gth}{N-2})$ if $N>2$ and $1\leq p<\infty$ if $N=2$.

Let us first assume that $g$ is bounded. Then $\abs{G(v)}\leq m\abs v$. 
Since $g$ is continuous, $G(v)\in L^1_\gs(\Gw)$ for any $v\in W^{1,2}_0(\Gw)$ and the functional $J$ is well defined and is of class $C^1$ in $W^{1,2}_0(\Gw)$. Furthermore
\bel{Y12}
\lim_{\norm v_{W^{1,2}}\to\infty}J(v)=+\infty.
\ee
Let $\{u_n\}$ be a minimizing sequence. 
By $(\ref{Y12})$, $\{u_n\}$ is bounded in $W^{1,2}_0(\Gw)$ and thus relatively compact in $L^1_\gs(\Gw)$ and in $L^2(\Gw)$. 
Hence  there exist $u\in L^2(\Gw)$ and $v\in L^1_\gs(\Gw)$ such that, up to a subsequence,  
$u_n\to v$ in $L^1_\gs(\Gw)$, and 
$u_n\to u$ strongly in $L^2(\Gw)$ and weakly in $W^{1,2}_0(\Gw)$. 
We can also assume that $u_n\to u$ $c_{1,2}$-quasi almost everywhere in the sense that there exists $E\subset \Omega$ with $c_{1,2}(E)=0$ such that $u_n(x)\to u(x)$ for any $x\in\Omega\backslash E$. 
According to \rprop{AdamsMazya2}, $\sigma$ is absolutely continuous with respect to the $c_{1,2}$-capacity. It follows that $\sigma(E)=0$ so that $u_n\to u$ $\sigma$-almost everywhere and thus $u=v$ $\sigma$-almost everywhere. Thus we have that $u_n\to u$ in $L^2(\Omega)$, in $L^1_\sigma(\Omega)$, $\sigma$-almost everywhere and weakly in  $W^{1,2}_0(\Gw)$. 
Then we have that  $\langle\gm, u_n\rangle\to \langle\gm, u\rangle$. By the dominated convergence theorem we have also that $G(u_n)\to G(u)$ in $L^1_\gs(\Gw)$.  Therefore 
\bel{Y13}\displaystyle J(u)\leq \liminf_{n\to\infty}J(u_n),\ee
 which implies that $u$ is a minimizer of $J$ in $W^{1,2}_0(\Gw)$.

If $g$ is unbounded, we write $g=g_1+g_2$ where $g_1=g\chi_{(-r_0,r_0)}$, $g_2=g\chi_{(-\infty-r_0]\cup [r_0,\infty)}$, where $r_0$ is defined in $(\ref{Z0})$. 
Hence 
$G(r)=G_1(r)+G_2(r)$ where $\abs{G_1(r)}\leq m\abs r$ and $G_2(r)$ is nonnegative. Using again $(\ref{Y4})$ we obtain that $(\ref{Y12})$ holds. A minimizing sequence $\{u_n\}$ inherits the same property as above, hence $u_n\to u$ $\gs$-almost everywhere in $\Gw$ and in $L^1_\gs(\Gw)$, this implies that
$G_1(u_n)\to G_1(u)$ in $L^1_\gs(\Gw)$ and $G_2(u)$ is $\gs$-measurable. By Fatou's lemma
$$\displaystyle\myint{}{}G_2(u)d\gs\leq \liminf_{n\to\infty}\myint{}{}G_2(u_n)d\gs,
$$
which implies that $(\ref{Y13})$ holds. Notice that, among the consequences, $X_G$ is closed subset of $W^{1,2}_0(\Gw)$. 
Hence $u$ in a minimizer of $J$ in $X_G(\Gw)$.\\
Uniqueness holds if $g$ is nondecreasing since it implies that $J$ is stricly convex and actually $X_G$ is a closed convex set.\smallskip

\nind {\it Step 2: The minimizer is a weak solution}. For $k>r_0$ we define $g_k$ by
$$g_k(r)=\left\{\BA{lll}g(r)\qquad&\text{if }\,\abs r\leq k\\
g(k)\qquad&\text{if }\,r> k\\
g(-k)\qquad&\text{if }\,r<-k
\EA\right.
$$
Then $g_k$ is continuous and bounded and the minimizer $u_k\in W^{1,2}_0(\Gw)$ of 
$$J_k(v)=\myfrac{1}{2}\myint{\Gw}{}\abs{\nabla v}^2dx+\myint{\Gw}{}G_k(v)\,d\gs-\langle\gm, v\rangle\,\;\text{ where }\; G_k(r)=\myint{0}{s}g_k(s) ds,
$$
is a weak solution (i.e. in the sense given by $(\ref{Y10-1})$) of 
\bel{Z1-k}\BA {lll}
-\Gd u+g_k(u)\gs=\gm\qquad\text{in }\Gw\\
\phantom{-\Gd +g_k(u)\gs}
u=0\qquad\text{on }\prt\Gw.
\EA
\ee
The following energy estimate holds
\bel{Z2-k}
\myint{\Gw}{}|\nabla u_k|^2dx+\myint{\Gw}{}u_kg_k(u_k)d\gs=\langle\gm,u_k\rangle\leq \norm\gm_{W^{-1,2}}\norm {u_k}_{W^{1,2}},
\ee
and it implies
\bel{Z3-k}
\myint{\Gw}{}|\nabla u_k|^2dx+\myint{\Gw}{}\abs {u_kg_k(u_k)}d\gs\leq \norm\gm^2_{W^{-1,2}}+m\gs(\Gw)=M,
\ee
for some $m=m(r_0)>0$. Up to a subsequence, $\{u_{k}\}_k$ converges to some $u$ as $k\to\infty$, weakly in  $W^{1,2}_0(\Gw)$, strongly in $L^{2}(\Gw)$, and 
almost everywhere in $\Gw$. By \rprop{AdamsMazya} the imbedding of $W^{1,2}(\Gw)$ in $L^q_\gs(\Gw)$ is compact for any $q<\frac{2\gth}{N-2}$. Hence the subsequence can be taken such that $u_{k}\to u$, $\gs$-almost everywhere as $k\to\infty$, and consequently 
$g_k(u_{k})\to g(u)$ $\gs$-almost everywhere. Let $E\subset \Gw$ be a Borel set, then for any $\gl>r_0$, 
$$\BA {lll}M\geq \myint{E}{}\abs{g_k(u_{k})u_{k}}d\gs\\[4mm]
\phantom{M}= \myint{E\cap\{|u_{k}|>\gl\}}{}\abs{g_k(u_{k})u_{k}}d\gs+\myint{E\cap\{|u_{k}|\leq \gl\}}{}\abs{g_k(u_{k})u_{k}}d\gs
\\[4mm]
\phantom{M}\geq \gl\myint{E\cap\{|u_{k}|> \gl\}}{}\abs{g_k(u_{k})}d\gs+\myint{E\cap\{|u_{k}|\leq \gl\}}{}\abs{g_k(u_{k})u_{k}}d\gs.
\EA$$
Therefore
$$\BA {lll}
\myint{E}{}\abs{g_k(u_{k})}d\gs=\myint{E\cap\{|u_{k}|>\gl\}}{}\abs{g_k(u_{k})}d\gs+\myint{E\cap\{|u_{k}|\leq \gl\}}{}\abs{g_k(u_{k})}d\gs
\\[4mm]\phantom{\myint{E}{}\abs{g_k(u_{k})}d\gs}
\leq \myfrac{M}{\gl}+\max\{\abs{g(r)}:\, |r|\leq \gl\}\gs(E)
\EA$$
For $\ge>0$ we first choose $\gl$ such that $\frac{M}{\gl}\leq\frac\ge 2$ and then $\gs(E)\leq \frac{\ge}{1+2\max\{\abs{g(r)}\leq \gl\}}$. This implies the uniform integrability of $\{g_k(u_{k})\}_k$ in $L^1_\gs(\Gw)$. Hence $g_k(u_{k})\to g(u)$ in $L^1_\gs(\Gw)$ by Vitali's convergence theorem. 
Since $u_{k}$ is a weak solution of $(\ref{Z1-k})$, there holds for any $\gz\in C^\infty_0(\Gw)$, 
\bel{Z4-k}
\myint{\Gw}{}\nabla u_{k}.\nabla\gz dx+\myint{\Gw}{}g_k(u_{k})\gz d\gs=\langle\gm,\gz\rangle.
\ee
Letting $k\to\infty$ we obtain, using the above convergence results, 
\bel{Z5-k}
-\myint{\Gw}{}\nabla u.\nabla\gz dx+\myint{\Gw}{}g(u)\gz d\gs=\langle\gm,\gz\rangle.
\ee
Hence $u$ is a weak solution. If $g$ is monotone, uniqueness is also a consequence of the weak formulation. Furthermore if $\gm,\gm'$ belong to 
$W^{-1,2}(\Gw)$ are such that $\gm-\gm'$ is a nonnegative measure, then $\langle \gm'-\gm,(u_\gm'-u_\gm)_+\rangle\leq 0$. Taking 
$(u_\gm'-u_\gm)_+$ for test function in the weak formulation yields $(u_\gm'-u_\gm)_+=0$.\qeda

 \subsection{The $L^1$ case}
 In the sequel we set
  \bel{X0}
  \BBX(\Gw)=\{\gz\in C^1(\overline\Gw), \gz=0\text{ on }\prt\Gw\text{ and }\Gd\gz\in L^\infty(\Gw)\},
\ee 
and $\BBX_+(\Gw)=\BBX(\Gw)\cap\{\gz\in C^1(\overline\Gw):\gz\geq 0\text{ in }\overline\Gw\}$.
We recall (see e.g. \cite{Ve1}) that if $f\in L_\gr^1(\Gw)$ and $u\in L^1(\Gw)$ is a very weak solution of 
 \bel{X1}-\Gd u=f\qquad\text{in }\,\Gw,
\ee 
there holds
\bel{X2}
-\myint{\Gw}{}\abs u\Gd\gz dx\leq \myint{\Gw}{} f\text{sign} (u)\gz dx\qquad\text{for all }\gz\in \BBX_+(\Gw),
\ee
and
\bel{X3}
-\myint{\Gw}{}u^+\Gd\gz dx\leq \myint{\Gw}{} f\text{sign}_+ (u)\gz dx\qquad\text{for all }\gz\in \BBX_+(\Gw).
\ee


\bprop{Bre} 
Assume $N\geq 2$, $\gs\in \CM^+^+_{\frac{N}{N-\gth}}(\Omega)$ with $N\geq \gth>N-2$ and $g:\BBR\mapsto\BBR$ is  a continuous nondecreasing function vanishing at $0$. 
If $\gm\in L^1_\gr(\Gw)$ there exists a unique $u:=u_\gm\in L^1(\Gw)$ very weak solution  of $(\ref{Z1})$. 
Furthermore, if $u_\gm,u_{\gm'}\in L^1(\Gw)$ are the very weak solutions of $(\ref{Z1})$ with right-hand sides $\gm,\gm'\in L^1_\gr(\Gw)$, then 
\bel{X4}\BA {lll}
-\myint{\Gw}{}\abs{u_\gm-u_{\gm'}}\Gd\gz dx+\myint{\Gw}{}\abs{g(u_\gm)-g(u_{\gm'})}\gz d\gs\leq \myint{\Gw}{} (\gm-\gm')\text{sign} (u_\gm-u_{\gm'})\gz dx,
\EA\ee
and
\bel{X5}\BA {lll}
-\myint{\Gw}{}(u_\gm-u_{\gm'})_+\Gd\gz dx+\myint{\Gw}{}(g(u_\gm)-g(u_{\gm'}))_+\gz d\gs\leq \myint{\Gw}{} (\gm-\gm')\text{sign}_+ (u_\gm-u_{\gm'})\gz dx
\EA\ee
for any $\gz\in \BBX_+(\Gw)$. 
In particular the mapping $\mu\to u_\mu$ is nondecreasing. 
\es

The following result will be used several time in the sequel. its proof is standard but we present it for the sake of completeness.

\blemma{quasi} Assume $N>q\geq1$ and $\gs\in \CM^+^+_{\frac{N}{N-\gth}}$ with $N\geq \gth>N-q$. Then $\gs$ vanishes on any Borel set with $c_{1,q}$-capacity zero.
\es
\Proof 
It suffices to prove the result when $E$ is compact. We define the $\Gl_\gth$ Hausdorff measure of a set $E$ by
\bel{AA1}\BA {lll} \displaystyle
\Gl_\gth(E)=\lim_{\gk\to 0}\Gl^\gk_\gth(E)
:=\lim_{\gk\to 0}\, \inf \left\{\sum_{j=1}^\infty r_j^\gth:\,0<r_j\leq\gk\leq\infty,\,E\subset\bigcup_{j=1}^\infty B_{r_j}(a_j)\right\}.
\EA\ee
Note that $\Gl^\infty_\gth(E)$ is the Hausdorff content of $E$ and it is smaller than $($diam$\,(E))^\gth$. 
For any covering of $E$ by balls $B_{r_j}(a_j)$, $j\geq 1$, we have 
$$ \gs(E) \leq \sum_{j=1}^\infty \sigma(B_{r_j}(a_j))
\le \|\sigma\|_{{\frac{N}{N-\gth}}} \sum_{j=1}^\infty  r_j^\theta. $$ 
It follows that 
$$ \gs(E) \leq \|\sigma\|_{\frac{N}{N-\gth}} \Gl_\gth(E).
$$
Next, if $c_{1,q}(E)=0$ then $\Gl_\gth(E)=0$ according to  \cite[Th. 5.1.13]{AdHe}, and thus $\sigma(E)=0$ by the previous inequality. \qeda\medskip

We introduce the flow coordinates near $\prt\Gw$ defined by 
$$\Gp(x)=\left(\gr(x),\gt(x)\right)\in [0,\ge_0]\ti\prt\Gw\quad\text{where }\,\gt(x)=proj_{\prt\Gw}(x).
$$
It is well-known that for $\ge_0$ small enough, $\Gp$ is a $C^1$-diffeomorphism from $\Gw_{\ge_0}:=\{x\in\overline\Gw:\gr(x)\leq\ge_0\}$ to 
$[0,\ge_0]\ti\prt\Gw$. With this diffeomorphism we can assimilate the surface measure $dS_\ge$ on $\Gs_\ge=\{x\in\Gw:\gr(x)=\ge\}$
with the surface measure $dS$ on $\Gs_0=\prt\Gw$ by setting
$$\myint{\Gs_\ge}{} v(x)dS_\ge(x)=\myint{\Gs_0}{} v(\ge,\gt)dS(\gt). 
$$

\blemma{tr} Assume $N\geq 2$ and $\gl\in \mathfrak M(\Gw)$ satisfies 
\bel{alpha-1}
\myint{\Gw}{}\gr d\abs\gm<\infty.
\ee
Then $u=\BBG[\gm]$ satisfies 
\bel{alpha-2}
\lim_{\ge\to 0}\myint{\Gs_0}{} |u|(\ge,\gt)dS(\gt)=0.
\ee
\es
\Proof If $u=\BBG[\gm]$, it is the unique weak solution of $-\Gd u=\gm$ in $\Gw$, $u=0$ on $\prt\Gw$. 
Hence 
$u=u_1-u_2$ where $u_1=\BBG[\gm^+]$ and $u_2=\BBG[\gm^-]$. Because $\gm_+$ and $\gm_-$ satisfies the integrability condition 
$(\ref{alpha-1})$ both $u_1$ and $u_2$ has a zero measure boundary trace ( $M$- boundary trace in the sense of \cite[Sec 1.3]{MaVe3}). Hence, taking for test function the function $\gz=1$, 
\bel{alpha-3}\displaystyle
\lim_{\ge\to 0}\myint{\Gs_0}{} u_j(\ge,\gt)dS(\gt)=0,
\ee
which implies $(\ref{alpha-1})$.
\qeda\medskip

This result allows us to obtain the uniqueness of the solution even if the right-hand side is a measure. 
\blemma{Uniq} Assume $N\geq 2$, $\gs\in \CM^+^+_{\frac{N}{N-\gth}}(\Omega)$ with $N\geq \gth>N-2$ and $g:\BBR\mapsto\BBR$ is  a continuous nondecreasing function. 
If $\gm\in \mathfrak M(\Gw)$ there exists at most one very weak solution  of $(\ref{Z1})$. 
\es
\Proof By \rlemma{quasi} with $\ga=1$, $p=2$,   $\gs$ is absolutely continuous with respect to the  $c_{1,2}$ capacity (it is diffuse in the terminology of \cite {BMP}), and if $h\in L^1_\gs(\Gw)$ the measure 
$h_+\gs$, which is the increasing limit of $\inf\{n,h_+\}\gs$ is also diffuse. Similarly $h_-\gs$ is diffuse and so is $h\gs$. Next we assume that 
$u$ and $u'$ are two very weak solutions of $(\ref{Z1})$ and set $w=u-u'$. Hence 
$$-\Gd w+(g(u)-g(u'))\gs=0.
$$
Since $\gr(g(u)-g(u'))\in L^1_\gs(\Gw)$, it follows from \rlemma {tr} that 
$$\displaystyle \lim_{\ge\to 0}\myint{\Gs_\ge}{}\abs w(\ge,\gt) dS(\gt)=0
$$
We use Kato inequality for measure as in 
\cite [Th 1.1]{BrPo}: Since $w\in L^1(\Gw)$, $\Gd w^+$ is a diffuse measure and
$$\Gd w^+\geq \chi_{_{\{w\geq 0\}}}\Gd w=\chi_{_{\{w\geq 0\}}}(g(u)-g(u'))\gs\geq 0 \;\text{ in }\Gw
$$
Since $w^+$ has a M-boundary trace by \rlemma {tr}, we can apply \cite[Lemmma 1.5.8]{MaVe3} with $\gm=-\chi_{_{\{w\geq 0\}}}(g(u)-g(u'))\gs$ which is a measure in $\mathfrak M_\gr(\Gw):=\{\gn\in \mathfrak M(\Gw):\gr\gn\in \mathfrak M_b(\Gw)\}$. Then there exists $\gt\in \mathfrak M^+_\gr(\Gw)$ such that 
$$-\Gd w^+=\gm-\gt.
$$
Equivalently
$$-\Gd w^++\chi_{_{\{w\geq 0\}}}(g(u)-g(u'))\gs=-\gt. 
$$
 Since the M-boundary trace of $w^+$ is zero, it follows that $w^+=-\BBG[\chi_{_{\{w\geq 0\}}}(g(u)-g(u'))\gs+\gt]$. Hence $w^+=0$ and 
 $u\leq u'$. Similarly $u'\leq u$.\qeda\medskip

The following variant will be useful in the sequel.

\blemma{Order} Assume $N\geq 2$, $\gs\in \CM^+^+_{\frac{N}{N-\gth}}(\Omega)$ with $N\geq \gth>N-2$ and $g:\BBR\mapsto\BBR$ is  a continuous nondecreasing function. If $u,u'\in L^1(\Gw)$ are such that $\gr g(u)$ and $\gr g(u')$ belong to $L^1_\gs(\Gw)$ and satisfy
\bel{AB0}
-\myint{\Gw}{}(u-u')\Gd\gz dx +\myint{\Gw}{}(g(u)-g(u'))\gz d\gs=\myint{\Gw}{}\gz d\gn\quad\text{for all }\,\gz\in \BBX_+(\Gw)
\ee
for some $\gn\in \mathfrak M_+(\Gw)$ diffuse with respect to the $c_{1,2}$-capacity, then $u\geq u'$ $c_{1,2}$-quasi everywhere in $\Gw$.
\es
\Proof We use Kato's inequality, \rlemma {tr} and \cite[Lemmma 1.5.8]{MaVe3} in the same way as in the proof of \rlemma {Uniq} since the measures $(g(u)-g(u'))d\gs$ and $\gn$ are diffuse,
$\Gd (u'-u)$ is diffuse, hence
$$\Gd (u'-u)_+\geq\chi_{_{\{u'\geq u\}}}\Gd (u'-u)=(g(')-g(u))\chi_{_{\{u'\geq u\}}}+\chi_{_\{u'\geq u\}}\gn\geq 0
$$
Since $u'-u\in W^{1,q}_0(\Gw)$ for any $1<q<\frac{N}{N-1}$, we conclude that $(u'-u)_+=0$ almost everywhere and $c_{1,2}$-quasi everywhere by \cite[Th 6.1.4]{AdHe}.\qeda
\medskip

The next result and the corollary which follows are the key-stone for the proof of \rprop{Bre}.

\blemma{Bre-2} Let $\gs\in \CM^+^+_{\frac{N}{N-\gth}}(\Omega)$ with $N\geq \gth>N-2$, $h\in L^\infty_\gs(\Gw)$, $f\in L^s(\Gw)$ with $s>\frac N2$ and $w\in L^1(\Gw)$ be the very weak solution of 
\bel{AB1}\BA {lll} \displaystyle
-\Gd w+h\gs=f\qquad&\text{in }\;\Gw\\
\phantom{-\Gd+h\gs}
w=0&\text{in }\;\prt\Gw.
\EA\ee
Then  $w$ is continuous in $\overline\Gw$ and for any nondecreasing bounded function $\gg\in C^{2}(\BBR)$ vanishing at $0$, 
 there holds
\bel{AB2}\BA {lll} \displaystyle
-\myint{\Gw}{}j(w)\Gd\gz dx+\myint{\Gw}{}\gg(w)h\gz d\gs\leq \myint{\Gw}{}\gg(w)\gz f dx\qquad\text{for all }\gz\in\BBX_+(\Gw),
\EA\ee
where $j(r)=\myint{0}{r}\gg(s) ds$.
\es

\noindent\Proof The solution is unique and expressed by $w=\BBG[f-h\gs]$. Since $\frac{N}{N-\theta}>\frac N2$, $w\in C^\ga(\overline\Gw)$ for some $\ga\in (0,1)$ by \rlemma{Miy}. Hence $\gg(w)$ is continuous and therefore measurable. We extend $\gs$ by zero in $\Gw^c$ and denote $\gs_n=\gs\ast\eta_n$ where $\{\eta_n\}$ is a sequence of mollifiers. Then 
$\gs_n\to\gs$ in the narrow topology of $\Gw$. For $n\in\BBN^*$, let $w_n$ be the solution of 
\bel{AB1-n}\BA {lll} \displaystyle
-\Gd w_n+h\gs_n=T_n(f)\qquad&\text{in }\;\Gw\\
\phantom{-\Gd_n+h\gs_n}
w_n=0&\text{in }\;\prt\Gw,
\EA\ee
where $T_n(f)=\min\{|f|,n\}$sgn$(f)$. Then $w_n\in W^{2,s}(\Gw)\cap W^{1,\infty}_0(\Gw)$ for all $1<s<\infty$. By Green's formula
\bel{AB2-n}\displaystyle
-\myint{\Gw}{}j(w_n)\Gd\gz dx+\myint{\Gw}{}\gg(w_n)h\gz d\gs\leq \myint{\Gw}{}\gg(w_n)\gz f dx\qquad\text{for all }\gz\in\BBX_+(\Gw).
\ee
Since $w_n\to w$ uniformly in $\overline\Gw$, $(\ref{AB2})$ follows.\qeda

\bcor{Bre-3} Under the assumptions of \rlemma{Bre-2}, there holds
\bel{AB5}\BA {lll} \displaystyle
-\myint{\Gw}{}\abs w\Gd\gz dx+\myint{\Gw}{}sign_0 (w)h\gz d\gs
\leq \myint{\Gw}{}sign_0 (w)\gz f dx,
\EA\ee
and
\bel{AB6}\BA {lll} \displaystyle
-\myint{\Gw}{}w_+\Gd\gz dx+\myint{\Gw}{}sign_+ (w)\gz hd\gs
\leq \myint{\Gw}{}sign_+ (w)\gz f dx, 
\EA\ee
for any $\gz\in\BBX_+(\Gw)$. 
Moreover there exists a constant $C>0$ depending only on $\Omega$ such that 
\bel{AB61}\BA {lll} \displaystyle
\myint{\Gw}{} sign_0 (w) h d\gs \leq C \myint{\Gw}{}|f| dx.  
\EA\ee 
\es

\noindent\Proof For proving $(\ref{AB5})$ we consider a sequence $\{\gg_k\}$ of odd nondecreasing functions such that 
$$\gg_k(r)=\left\{\BA {lll}
\phantom{-}1\quad&\text{if }\;r\geq 2k^{-1}\\
\phantom{-}0\quad &\text{if }\;-k^{-1}\leq r\leq k^{-1}\\
-1\quad&\text{if }\;r\leq -2k^{-1}
\EA\right.
$$
and such that $\{r\gg_k(r)\}$ is nondecreasing for any $r$. 
Using $\gamma_k$ in place of $\gamma$ in $(\ref{AB2})$ we obtain 
\begin{equation}\label{AB20} 
  -\myint{\Gw}{}j_k(w)\Gd\gz dx+\myint{\Gw}{}\gg_k(w)\gz  hd\gs
\leq \myint{\Gw}{}\gg_k(w)\gz f dx\qquad\text{for all }\gz\in\BBX_+(\Gw),
\end{equation} 
where $j_k(r)=\myint{0}{r}\gg_k(s) ds$. Since $\gg_k(w)\uparrow w$ on $\Gw_+:=\{x\in\Gw:w(x)>0\}$, 
there holds by the monotone convergence theorem,
$$
\myint{\Gw_+}{}\gg_k (w)\gz \abs hd\gs\uparrow \myint{\Gw_+}{}w\gz\abs h d\gs\qquad\text{as }\,k\to\infty.
$$
Since 
$$\abs{\myint{\Gw_+}{}(w-\gg_k(w))\gz hd\gs}\leq \myint{\Gw_+}{}\abs{(w-\gg_k(w))\gz h}d\gs=\myint{\Gw_+}{}(w-\gg_k(w))\gz |h|d\gs,
$$
we obtain
$$
\myint{\Gw_+}{}\gg_k (w)h\gz d\gs\to \myint{\Gw_+}{}wh\gz d\gs\qquad\text{as }\,k\to\infty.
$$
Similarly, $\gg_k(w)\downarrow w$ on $\Gw_-:=\{x\in\Gw:w(x)<0\}$ so that 
$$
\myint{\Gw_-}{}\gg_k (w)h\gz d\gs\to \myint{\Gw_-}{}wh\gz d\gs\qquad\text{as }\,k\to\infty.
$$
Combining these two results yields
$$\myint{\Gw}{}\gg_k (w)\gz hd\gs\to \myint{\Gw_+}{}w\gz hd\gs-\myint{\Gw_-}{}w\gz hd\gs=\myint{\Gw}{}sgn_0(w)\gz hd\gs.
$$
Usiing dominated convergence theorem there holds
$$\myint{\Gw}{}\gg_k (w)\Gd \gz dx 
\to \myint{\Gw}{}sgn_0(w)\Gd \gz dx, 
$$
and 
$$\myint{\Gw}{}\gg_k (w) \gz f dx\to \myint{\Gw}{}sgn_0(w)\gz f dx.
$$
This implies $(\ref{AB5})$. The proof of $(\ref{X4})$ is similar.

Eventually we prove $(\ref{AB61})$. Let $\eta_1$ be the solution of 
\bel{A16}\BA {lll} -\Gd \eta_1= 1 \qq &\text{in } \Gw\\
\phantom{-\Gd  }
\eta_1=0\qq &\text{in } \prt\Gw.
\EA\ee
Then $\eta_1=\BBG[1]\in \BBX_+(\Gw)$ and there exists $c,c'>0$ depending only on $\Omega$ such that $c\rho\leq\eta_1\leq c'\rho$. 
Given $\alpha\in (0,1]$, let $j_\ge(r)=(r+\ge)^\alpha-\ge^\alpha$, $r\geq 0$, and $\zeta=j_\ge(\eta_1)$. Note that $\zeta\in C^2(\overline{\Omega})$, 
$0\leq \zeta\le \eta^\alpha$, $\zeta=0$ on $\partial\Omega$, $j_\ge'>0$, $j_\ge''<0$, so that $-\Delta\zeta = j_\ge'(\eta_1)-j_\ge''(\eta_1)|\nabla\eta_1|^2 \geq 0$. We deduce from $(\ref{AB5})$ that 
$$ \int_\Omega sgn_0(w) (\eta+\ge)^\alpha hd\sigma 
\leq \int_\Omega sgn_0(w)\eta^\alpha |f| dx 
+ \ge^\alpha\int_\Omega sgn_0(w) hd\sigma. $$ 
We obtain
$$ \int_\Omega sgn_0(w) \rho^\alpha hd\sigma 
\leq C \int_\Omega \rho^\alpha |f| dx + \ge^\alpha |\tilde\sigma(\Omega)| $$ 
Letting $\ge\to 0$ and then $\alpha\to 0$ we infer the result by dominated convergence. 

\qeda

\medskip

We are now in position to prove \rprop{Bre}. 

\medskip

\nind{\it Proof of \rprop{Bre}}. We divide the proof into several steps. 

\smallskip 

\nind  {\it Step 1: We assume that $\gm\in L^\infty(\Gw)$}.
Let $\{\eta_n\}$ be a sequence of molifiers and $\gs_{n}=\gs\ast\eta_n$. If $\gm\in L^\infty(\Gw)$, the solution $u_n=u_{n,\gm}$ of
\bel{A14}\BA {lll} -\Gd u_n + g(u_n)\gs_{n} = \gm \qq &\text{in } \Gw\\
\phantom{-\Gd  + g(u_n)\gs_{n} }
u_n=0\qq &\text{in } \prt\Gw,
\EA\ee
is continuous in $\overline\Gw$.  Since 
\bel{A14-1}
-\BBG[\gm^-]\leq -u^-_n\leq 0\leq u^+_n\leq \BBG[\gm^+]
\ee
by the maximum principle, the sequence $\{u_n\}$ is uniformly bounded. 
Recalling that  $g$ is nondecreasing we have that the sequence $\{g(u_n)\}$ is also 
uniformly bounded in $\Gw$, hence $g(u_n)\gs_{n}$ is bounded in 
$\CM^+_{\frac{N}{N-\gth}}(\Gw)$ independently of $n$, and from $(\ref{Y2+3})$ it implies that $u_n$ is bounded in $C^\ga(\overline\Gw)$ for some $\ga\in (0,1]$ independently of $n$. Up to some subsequence, $\{u_n\}$, and thus also $\{g(u_n)\}$, are then uniformly convergent in $\overline\Gw$ with limit $u=u_{\gm}$ and $g(u)=g(u_\gm)$. Because $\gs\ast\eta_n$ converges to $\gs$ in the narrow topology, $u_{\gm}$ is a very weak solution of  $(\ref{Z1})$. Notice that being continuous, $g(u)$ is measurable for the measure $\gs$. By \rlemma{Uniq}, $u_{\gm}$ is the unique solution of $(\ref{Z1})$, hence the whole sequence $\{u_{\gm_n}\}$ converges to $u_{\gm}$.
Applying  \rcor{Bre-3} with $w=u$, $\tilde\gs=\gs$ and $\zeta=\eta_1$ yields 
\bel{A14'}\BA {lll} \displaystyle
\myint{\Gw}{}\abs {u} dx+\myint{\Gw}{}\abs{g(u)}\eta_1 d\gs\leq \myint{\Gw}{}\abs{\gm} \eta_1dx,
\EA\ee
and $(\ref{AB6})$ with $\zeta=\eta_1$ gives  �
\bel{A18}\BA {lll} \displaystyle
\myint{\Gw}{}(u-u')_+ dx +\myint{\Gw}{}(g(u)-g(u'))_+\eta_1 d\gs
\leq \myint{\Gw}{}\eta_1 sign_+ (u-u') (\gm-\gm')_+ dx.
\EA\ee
which implies the monotonicity of the mapping $\gm\mapsto u_\gm$.\smallskip

\nind{\it Step 2: We assume that $\gm\in L^1(\Gw)$ is bounded from below}.
Set $\ell=\rm{ess}\,\inf\gm$. For $k>0$ set $\gm_{k}=\min\{k,\gm\}$ and $u_{k}:=u_{\gm_{k}}\in L^\infty(\Omega)$. 
The sequence $\{\gm_{k}\}$ is nondecreasing, hence
according to Step 1, the sequence $\{u_{k}\}$ is a nondecreasing sequence of continuous functions in $\overline\Gw$ bounded from below by $\ell\eta_1$, where $\eta_1$ is defined in $(\ref{A16})$. 
Its pointwise limit, denoted by $u$ is thus lower semicontinuous. 
Moreover $g(u_{k})\to g(u)$ pointwise, hence $g(u)$ is lower semicontinuous and thus $\sigma$-measurable. 
Relation $(\ref{A14'})$ applied to $\mu_k$ and $u_k$ gives 
$$ \myint{\Gw}{}\abs {u_k} dx+\myint{\Gw}{}\abs{g(u_k)}\eta_1 d\gs
\leq \myint{\Gw}{}\abs{\gm_k} \eta_1dx.$$ 
Passing to the limit using Fatou's lemma in the  left-hand side and the dominated convergence theorem in the right-hand side yields 
\begin{equation}\label{Eq30}
  \myint{\Gw}{}\abs {u} dx+\myint{\Gw}{}\abs{g(u)}\eta_1 d\gs
\leq \myint{\Gw}{}\abs{\gm} \eta_1dx. 
\end{equation} 
We deduce that $u\in L^1(\Gw)$ and $\gr g(u)\in L^1_\sigma(\Gw)$. We have indeed a more precise result. Since $g$ vanishes at $0$ 
$g(u_k)=g(u^+_k)+g(-u^-_k)$. Hence $\gr g(u^+_k)\to \gr g(u^+)$ in $L^1_\sigma(\Gw)$ by the monotone convergence theorem. Furthermore 
$g(-u^-_1)\leq g(-u^-_k)\leq 0$, which implies that $\gr g(-u^-_k)\to \gr g(-u^-)$ in $L^1_\sigma(\Gw)$ by the dominated convergence theorem which finally implies that $\gr g(u_k)\to \gr g(u)$ in $L^1_\sigma(\Gw)$.
Using $\gz\in\BBX_+(\Gw)$ as a test function in the very weak formulation of the equation satisfied by $u_k$ gives 
$$ -\myint{\Gw}{}u_{k} \Delta\zeta dx+\myint{\Gw}{}g(u_{k})\zeta d\gs
= \myint{\Gw}{}\zeta\gm_{k} dx. $$
Since $u_k\to u$ almost everywhere and $-l\eta_1\leq u_k\leq u$ with $u\in L^1(\Omega)$,  we can pass to the limit in the first term to obtain 
$\myint{\Gw}{}u_{k} \Delta\zeta dx\to \myint{\Gw}{}u \Delta\zeta dx$. 
Because $|\mu_k|\leq |\mu|\in L^1(\Omega)$ and $\mu_k\to \mu$ almost everywhere,  we can also  pass 
to the limit in the last term: $\myint{\Gw}{}\zeta\gm_k dx\to \myint{\Gw}{}\zeta\gm dx$. Since $\gz g(u_k)\to \gz g(u)$ in $L^1_\sigma(\Gw)$ we conclude that 
It remains to pass to the limit in the nonlinearity. Because $u_k\uparrow u$ and $g$ is nondecreasing, we have $g(u_k)\uparrow g(u)$. Thus by the monotone convergence theorem, 
\begin{eqnarray*} 
-\myint{\Gw}{}u \Delta\zeta dx+\myint{\Gw}{}g(u)\zeta d\gs
= \myint{\Gw}{}\zeta\gm dx, \end{eqnarray*} 
and $u$ is very weak solution of $(\ref{Z1})$. 

\smallskip 

\nind{\it Step 3: We assume that $\gm\in L^1(\Gw)$}. 
For $\ell \in\BBR$, we set $\gm^{\ell}=\sup\{\gm,\ell\}$ and denote by $u^{\ell}$  the solution of $(\ref{Z1})$ with right-hand side $\gm^{\ell}$. 
Note that the sequence $\{\gm^{\ell}\}_\ell$ is increasing, bounded from above by $\mu^+$  so that 
$u^{\ell}\leq  u_{\gm^+}$, where 
$u_{\gm^+}$ is the solution of $(\ref{Z1})$ with right-hand side  $\gm^+$ which exists according to the previous step, and 
the sequence $\{u^{\ell}\}_\ell$ is monotone nondecreasing with $\ell$ with pointwise limit $u$ when $\ell\to-\infty$. Hence
 $u\leq u^{\ell}\leq  u_{\gm^+}$ for any $\ell \leq 0$. The sequence  $\{g(u^{\ell})\}_\ell$ is monotone nondecreasing with limit $g(u)$ when $\ell\to-\infty$, and there holds  $g(u)\leq g(u^{\ell})\leq  g(u_{\gm^+})$ for any $\ell \leq 0$. Since $g(u^{\ell})$ is lower semicontinuous and $\gs$-measurable, $g(u)$ shares the same properties.
 
\nind Applying $(\ref{Eq30})$ to $\mu=\mu^{\ell}$ and $u=u^{\ell}$ gives 
$$  \myint{\Gw}{}\abs {u^{\ell}} dx+\myint{\Gw}{}\abs{g(u^{\ell})}\eta_1 d\gs
\leq \myint{\Gw}{}\abs{\mu^{\ell}} \eta_1dx.  $$ 
Passing to the limit in the  left-hand side using Fatou's lemma we obtain 
$$  \myint{\Gw}{}\abs {u} dx+\myint{\Gw}{}\abs{g(u)}\eta_1 d\gs
\leq \myint{\Gw}{}\abs{\mu} \eta_1dx.  $$ 
We deduce that $u\in L^1(\Gw)$ and $\gr g(u)\in L^1_\sigma(\Gw)$. 
We conclude as in Step 2 that $u$ is solution of $(\ref{Z1})$. 

\smallskip

\nind{\it Step 4: Proof of $(\ref{X4})$ and $(\ref{X5})$.} 

For $\ell<0<k$ we set $\gm^\ell_k=\sup\{\ell,\inf\{k,\gm\}\}$ and denote by $u_k^{\ell}$ the solution of $(\ref{Z1})$ with right-hand side $\gm^\ell_k$. Then, by \rcor{Bre-3}, for any $\gz\in\BBX(\Gw)$ there holds
$$\BA {lll}
-\myint{\Gw}{}\abs{u_k^{\ell}-u_{k'}^{\ell'}}\Gd\gz dx+\myint{\Gw}{}\abs{g(u_k^{\ell})-g(u_{k'}^{\ell'})}\gz d\gs 
\leq \myint{\Gw}{}\rm {sign}_0(w)(\gm^\ell_k-\gm_{k'}^{\ell'})\gz dx.
\EA$$
Using the previous convergence theorem when $k\to\infty$ and then $\ell\to-\infty$, we derive $(\ref{X4})$. The proof of $(\ref{X5})$ is similar.
\qeda\medskip

\nind\Remark If it is not assumed that $g$ is nondecreasing, the above proof by monotonicity does not work. However the existence will follow from 
Theorem B if it is assumed that the extra assumptions in this theorem are satisfied: $\gth>N-q$ for some $q\in (1,\frac{N}{N-1})$ and the growth assumptions of Theorem B.

\subsection{Diffuse case}

We recall that a measure $\gm$ is said to be diffuse with respect to the $c_{s,p}$-capacity defined in $(\ref{DefCapBessel} )$ if $\abs\gm$ vanishes on all sets with zero $c_{s,p}$-capacity.
An important result due to Feyel and de la Pradelle \cite{FePr}  is the following: 

\bprop{FeyelPradelle} Let $\alpha > 0$ and $1 < p < \infty$. 
If $\lambda \in \mathfrak M^+_b(\Gw)$ does not charge sets with zero 
$c_{\alpha,p}$-capacity, there exists an increasing sequence 
$\{\lambda_n\}\subset  H^{-\alpha,p'}(\Gw)\cap\mathfrak M_b^+(\Gw)$, $\lambda_n$ with
compact support in $\Omega$, which converges to $\lambda$.
\es

\bprop{diffuse} Assume $\gs\in \CM^+^+_{\frac{N}{N-\gth}}$ with $N\geq \gth>N-2$, 
 and that $g:\BBR\mapsto\BBR$ is  a continuous nondecreasing function vanishing at $0$.
Then for any $\gm\in \mathfrak M^+_b(\Gw)$ diffuse with respect to the $c_{1,2}$-capacity there exists a unique very weak solution $u$ to $(\ref{Z1})$.
\es

\noindent\Proof 
According to \rprop{FeyelPradelle}, there exists  an increasing sequence of nonnegative measures $\{\gm_n\}$ belonging to $W^{-1,2}(\Gw)$ and converging to $\gm$ and by \rprop{varia}, $\{u_{\gm_n}\}$ is a nondecreasing sequence of weak solutions of $(\ref{Z1})$ with $\gm=\gm_n$.
We claim that $u_{\gm_n}\uparrow u_{\gm}$ which is a very weak solution of (\ref{Z1}).  
There holds,  
\begin{eqnarray*} 
\myint{\Gw}{}u_{\mu_n} dx+\myint{\Gw}{}g(u_{\mu_n})\eta_1d\gs
= \myint{\Gw}{}\eta_1d\gm_n 
\leq \myint{\Gw}{}\eta_1d\gm,
\end{eqnarray*} 
where $\eta_1$ is defined in $(\ref{A16})$. Since $u_{\gm_n}\geq 0$,  $u_{\gm_n}\uparrow u$ and $g(u_{\gm_n})\uparrow g(u)$. Since
$u_{\gm_n}$ is $\gs$-measurable by \rprop{varia}, $u$ is also $\gs$-measurable. Hence $g(u)$ shares this measurability property since  $g$ is continuous. Hence, by the monotone convergence theorem
 \bel{Z4-k'}
\myint{\Gw}{}u dx+\myint{\Gw}{}g(u)\eta_1d\gs
= \myint{\Gw}{}\eta_1d\gm.
\ee
Furthermore $u_{\gm_n}\to u$ in $L^1(\Gw)$.
 Indeed it suffices to show that $\{u_{\gm_n}\}$ is uniformly equiintegrable which follows from $0\leq \int_\omega u_{\gm_n}dx \le \int_\omega u dx$ and the fact that $u\in L^1(\Omega)$. We show in the same way that $\gr g(u_{\gm_n})\to \gr g(u)$ in $L^1_\gs(\Gw)$. This implies that $u=u_\gm$ is the very weak solution of $(\ref{Z1})$.\qeda

\subsection{Subcritical nonlinearities: proof of Theorem B.}

\blemma{meas1} Assume $N>2$ and $\gs\in \CM^+^+_{\frac{N}{N-\gth}}(\Gw)$ with 
$N\geq \gth>N-2$. If $\gm\in\frak M_b(\Gw)$ and $\gl\geq 0$,  we set
$E_\gl[\gm]:=\{x\in\Gw:\BBG[\abs\gm](x)>\gl\}$. Then 
\bel{B1}\BA {lll} \displaystyle
e_\gl^\gs(\gm):=\myint{E_\gl[\gm]}{}d\gs
\leq c\norm\gm_{\frak M_b}^{\frac{\gth}{N-2}}\gl^{-\frac{\gth}{N-2}}\qquad\text{for all }\gl>0.
\EA\ee
\es

\noindent\Proof 
It suffices to prove the result if $\mu\geq 0$. Indeed since 
$\BBG[\abs\gm] = \BBG[\gm^+] + \BBG[\gm^-]$, we have 
$ E_\gl[\gm] \subset E_{\gl/2}[\gm^+] \cup E_{\gl/2}[\gm^-]$ 
and thus 
$e_\gl^\gs(\gm) \leq e_{\gl/2}^\gs(\gm^+) +e_{\gl/2}^\gs(\gm^+)$. 
If the result holds for nonnegative measure, in particular for $\mu^\pm$, then 
\begin{eqnarray*} 
\gl^\frac{\gth}{N-2} e_\gl^\gs(\gm) 
& \leq & 
c (\mu^+(\Omega)^{\frac{\gth}{N-2}} + \mu^-(\Omega)^{\frac{\gth}{N-2}}) 
 \leq   
 c (\mu^+(\Omega) + \mu^-(\Omega))^\frac{\gth}{N-2} \\ 
& = & c\norm\gm_{\frak M_b}^{\frac{\gth}{N-2}}. 
\end{eqnarray*} 
Thus, we assume from now on that $\mu$ is nonnegative.  

If $\gm=\gd_a$ for some $a\in\Gw$, then 
$\BBG[\gd_a](x)\leq c_{_N}\abs{x-a}^{2-N}$ so that 
$E_\gl[\gd_a]\subset B_{(\frac{c_{_N}}{\gl})^{\frac{1}{N-2}}}(a)$. 
Since $\gs\in \CM^+^+_{\frac{N}{N-\gth}}(\Gw)$  it follows that 
\begin{equation}\label{Equ3} 
e_\gl^\gs(\gd_a)\leq c\gl^{-\frac{\gth}{N-2}}.
\end{equation} 
Let $E\subset\Gw$ be a Borel set. 
For any given $t>0$ there holds 
$$ \myint{E}{}\BBG[\gd_a]d\gs=\myint{E\cap E_t[\gd_a] }{}\BBG[\gd_a]d\gs
+\myint{E\cap E^c_t[\gd_a] }{}\BBG[\gd_a]d\gs. $$ 
Clearly $ \myint{E\cap E^c_t[\gd_a] }{}\BBG[\gd_a]d\gs \le t\sigma(E)$ and 
$$ \myint{E\cap E_t[\gd_a] }{}\BBG[\gd_a]d\gs
\le \myint{E_t[\gd_a] }{}\BBG[\gd_a]d\gs
\le -\myint{t}{\infty}s\,de_s^\gs(\gd_a)
\le c\myfrac{\gth t^{1-\frac{\gth}{N-2}}}{\gth+2-N}, $$ 
where the last inequality follows by integration by parts and the help of $(\ref{Equ3})$. 
Then 
$$ \myint{E}{}\BBG[\gd_a]d\gs
\leq t\sigma(E)+c\myfrac{\gth t^{1-\frac{\gth}{N-2}}}{\gth+2-N}. $$ 
Minimizing the  right-hand side with respect to $t$, we infer
\bel{B2}\BA {lll} \displaystyle
\myint{E}{}\BBG[\gd_a]d\gs\leq c\sigma(E)^{1-\frac{N-2}{\gth}}.
\EA\ee

We first suppose that $\gm=\sum_{j=1}^\infty\ga_j\gd_{a_j}$ for some $\ga_j>0$ and $a_j\in\Gw$. In particular $\sum_{j=1}^\infty\ga_j=\norm\gm_{\frak M^b}$. 
Using Fubini's theorem and $(\ref{B2})$ we see that for any Borel set $E\subset\Gw$, 
\bel{B3}\BA {lll}\displaystyle
\myint{E}{}\BBG[\gm](x)d\gs(x)
=\sum_{j=1}^\infty\ga_j\myint{E}{}\BBG[\gd_{a_j}(x)]d\gs(x)
\leq c\sigma(E)^{1-\frac{N-2}{\gth}}\norm\gm_{\frak M^b}.
\EA\ee
Taking in particular $E=E_\gl[\gm]$ we obtain 
$$\gl e_{\gl}^\gs(\gm)
\leq \myint{E_\gl[\gm]}{}\BBG[\gm](x)d\gs(x)
\leq c(e_{\gl}^\gs(\gm))^{1-\frac{N-2}{\theta}}\norm\gm_{\frak M^b},
$$
which implies the claim. 
Notice that the constant $c$ in the right-hand side depends only on $N$ and 
$\|\sigma\|_{\CM^+_{\frac{N}{N-\gth}}}$.

For  a general nonnegative measure $\gm\in\frak M_b(\Gw)$, we consider a sequence of nonnegative measures $\{\mu_n\}\subset \frak M_b(\Gw)$ where each $\mu_n$ is a sum of Dirac masses as before and such that $\mu_n\to \mu$ weakly as $n\to\infty$. 
Then we have 
$$ e_\gl^\gs(\gm_n):=\myint{E_\gl[\gm_n]}{}d\gs
\leq c\|\gm_n\|_{\frak M_b}^{\frac{\gth}{N-2}}\gl^{-\frac{\gth}{N-2}}, $$ 
with $\|\mu\|_{\frak M_b}\leq \displaystyle\liminf_{n\to\infty} \|\mu_n\|_{\frak M_b}$. 
We thus need to prove that 
\begin{equation}\label{Equ5} 
\liminf \myint{E_\gl[\gm_n]}{}d\gs \geq  \myint{E_\gl[\gm]}{}d\gs. 
\end{equation} 
We first observe that for any $t>0$ and $x\in \Omega$ the set 
$\{y\in\Omega:\, \BBG(x,y)>t\}$ is open (with $\BBG(x,x)=+\infty$). 
It follows from \cite{Bill}[Thm 2.1] that 
$\displaystyle\liminf_{n\to\infty} \mu_n(\{\BBG(x,\cdot)>t\})\geq \mu(\{\BBG(x,\cdot)>t\})$. 
We can take the $\liminf$ using Fatou's lemma in 
$$ \myint{\Gw}{} \BBG(x,y)\,d\mu_n(y) 
= \myint{0}{+\infty}  \mu_n(\{\BBG(x,\cdot)>t\})\,dt,  $$ 
to derive 
$$   \displaystyle\liminf_{n\to\infty}  \BBG[\mu_n](x) 
\geq \myint{0}{+\infty}  \mu(\{G(x,\cdot)>t\})\,dt 
= \myint{\Gw}{} G(x,y)\,d\mu(y) = \BBG[\mu](x). 
$$ 
We infer that for any $x\in\Omega$ such that $\chi_{E_\lambda(\mu)}(x)=1$ we have  $\displaystyle\liminf_{n\to\infty}\BBG[\mu_n](x) >\lambda$, 
hence $\BBG[\mu_n](x) >\lambda$ for $n$ large enough. 
Thus $\chi_{E_\lambda(\mu_n)}(x)=1$ eventually, and then 
$$ \displaystyle\liminf_{n\to\infty} \chi_{E_\lambda[\mu_n]}(x) \geq \chi_{E_\lambda[\mu]}(x) \qquad   
\text{for all } x\in\Omega. $$   
The claim $(\ref{Equ5})$ follows by Fatou's lemma. 
\qeda

\medskip

We are now in position to prove Theorem B. 

\medskip

\nind{\it Proof of Theorem B}. We note that if $g$ is nondecreasing, uniqueness follows from estimate \rlemma{Uniq}. 
Let $\{\eta_n\}$ be a sequence of mollifiers, $\gm_n=\gm\ast\eta_n$ and 
$u_n\in W_0^{1,2}(\Omega)$ a minimizing weak solution of 
\bel{B5}\BA {lll} -\Gd u_n + g(u_n)\gs = \gm_n \qq &\text{in } \Gw,\\
\phantom{-\Gd  + g(u_n)\gs}
u_n=0\qq &\text{in } \prt\Gw,
\EA\ee
given by \rprop{varia}. 
We write $g(r)=g_1(r)+g_2(r)$  with $g_1=g\chi_{(-r_0,r_0)}$, $g_2=g\chi_{(-\infty-r_0]\cup [r_0,\infty)}$, and set $m=\sup\{g(r):-r_0\leq r\leq r_0\}\geq 0$ and $m'=\inf\{g(r):-r_0\leq r\leq r_0\}\leq 0$. 
Then 
$$ -\BBG[\gm_n^-]-m\BBG[\sigma] \leq u_n\leq\BBG[\gm^+_n]-m'\BBG[\sigma]. $$ 
Since $\gs\in \CM^+^+_p(\Gw)$ for some $p>N/2$, 
$\BBG[\sigma]\in C^{0,\alpha}(\overline{\Omega})$ by  \rlemma{Miy}. 
Moreover $ \BBG[|\gm_n|]\in C(\overline{\Omega})$ since $|\mu_n|\in C(\overline{\Omega})$. 
It follows that 
\begin{equation}\label{eq40} 
    |u_n|\leq \BBG[|\gm_n|]+M \leq c_n, 
\end{equation} 
where $M,c_n\geq 0$. 

Since $u_n\in W^{1,2}_0(\Omega)$, its precise representative (that we identify with $u_n$) is
defined $c_{1,2}$-quasi-everywhere, is $c_{1,2}$-continuous and 
$$ u_n(x) = \lim_{r\to 0} \frac{1}{|B_r(x)|} \int_{B_r(x)} u_n(y)\,dy $$ 
for any $y\in \Omega\setminus E_n$ with $c_{1,2}(E_n)=0$ (see \cite{AdHe}). 
It follows that $|u_n|\le c_n$ in $E:=\cup E_n$. Note that $c_{1,2}(E)=0$ so that $\sigma(E)=0$ by \rlemma{quasi}. 
Hence $|u_n|\le c_n$ $\sigma$-almost everywhere, $g(u_n)\in L^\infty_\gs(\Omega)$,
and therefore $g(u_n)\sigma\in \CM^+^+_{\frac  N{N-\gth}}(\Gw)$. 
We can then apply \rcor{Bre-3} to obtain, for any $\gz\in\BBX_+(\Gw)$, that 
$$ 
-\myint{\Gw}{}\abs u_n\Gd\gz dx+\myint{\Gw}{}sign_0 (u_n)g(u_n)\gz d\gs
\leq \myint{\Gw}{}sign_0 (u_n)\gz \mu_n dx,
$$ 
which implies
\bel{B50}\BA {lll} \displaystyle
-\myint{\Gw}{}\abs u_n \Gd\gz dx+\myint{\Gw}{}|g_2(u_n)|\gz d\gs
\leq \myint{\Gw}{}sign_0 (u_n)\gz \mu_n dx + c\myint{\Gw}{}\gz d\gs.  
\EA\ee
We take $\zeta=\eta_1$ and obtain 
\bel{B7}\BA {lll} 
\myint{\Gw}{}\abs {u_n} dx+\myint{\Gw}{}\abs{g_2(u_n)}\eta_1d\gs
\leq \myint{\Gw}{}\abs{\gm_n}\eta_1dx+c \\[4mm]
\phantom{\myint{\Gw}{}\abs {u_n} dx+\myint{\Gw}{}\abs{g(u_n)}\eta_1d\gs}
\leq \myint{\Gw}{}\eta_1d\abs{\gm}+c = c',
\EA\ee
so that $\{u_n\}$ is bounded in $L^1(\Omega)$. 
We also have from \rcor{Bre-3} that 
$$ \myint{\Gw}{} sign_0 (u_n) g(u_n) d\gs \leq C \myint{\Gw}{}|\mu_n|\rho dx  $$ 
and so 
\begin{equation}\label{Equ51} 
  \myint{\Gw}{} |g_2(u_n)| d\gs
\leq C \myint{\Gw}{}|\mu_n| dx + \int_\Omega |g_1(u_n)| d\sigma \leq C
\end{equation} 
wth $C$ independent of $n$. 
We deduce that the sequence of measures $\{g(u_n)\}$ is bounded.

By the standard a regularity estimates, the sequence $\{u_n\}$ is bounded in  
$W^{1,q}(\Omega)$, $q<\frac{N}{N-1}$. 
Then there exists $u\in W^{1,q}(\Omega)$, $q<\frac{N}{N-1}$, such that, up to a subsequence, $u_n\to u$ in $L^1(\Omega)$ and also pointwise in $\Omega\setminus E$ where $c_{1,q}(E)=0$.  
We fix $q\in \left(1,\frac{N}{N-1}\right)$ such that $\gth>N-q$. 
In view of  \rlemma{quasi}, $\sigma(E)=0$ so that $g(u_n)\to g(u)$ $\sigma$-almost everywhere. 
Applying Fatou's lemma in $(\ref{Equ51})$ gives that $g(u)\in L^1_\sigma(\Omega)$.

In order to prove the uniform integrability of $\{g(u_n)\}$ for the measure $\gs$ we can assume that $\abs{g_2}\leq \tilde g$ with  a function satisfying $(\ref{Z7})$ still denoted by $\tilde g$ and let $E\subset \Gw$ be a Borel set. Then
$$\BA {lll}
\myint{E}{}\abs{g_2(u_n)}d\gs\leq \myint{E\cap\{\abs{u_n}\leq t\}}{}\abs{g_2(u_n)}d\gs+
\myint{E\cap\{\abs{u_n}> t\}}{}\abs{g_2(u_n)}d\gs
\\[4mm]
\phantom{\myint{E}{}\abs{g_2(u_n)}d\gs}
\leq \tilde g(t)\myint{E}{}d\gs+\myint{\{\abs{u_n}> t\}}{}\tilde g(\abs{u_n}) d\gs.
\EA$$
Then we estimate the second integral in the  right-hand side: for $\gl>M$ 
we set 
$$S_n(\gl)=\{x\in\Gw:\abs{u_n(x)}>\gl\}\quad\text{and }\;
b^\gs_n(\gl)=\myint{S_n(\gl)}{}d\gs.
$$ 
In view of $(\ref{eq40})$ we have $|u_n|\le \BBG(|\mu_n|) + M$
so that $S_n(\lambda)\subset E_{\lambda-M}[\mu_n]$. 
Hence $b^\gs_n(\gl)\leq e^\gs_{\gl-M}(\abs{\gm_n})$.  
This implies
$$\BA {lll}
\myint{\{\abs{u_n}> t\}}{}\tilde g(\abs{u_n}) d\gs
= -\myint{t}{\infty}\tilde g(\gl)d b^\gs_n(\gl)\\[4mm]
\phantom{\myint{\{\abs{u_n}> t\}}{}\tilde g(\abs{u_n}) d\gs}
\leq\myint{t}{\infty}b^\gs_n(\gl) d\tilde g(\gl)\\[4mm]
\phantom{\myint{\{\abs{u_n}> t\}}{}\tilde g(\abs{u_n}) d\gs}
\leq \myint{t}{\infty}e^\gs_{\gl-M}(\abs{\gm_n})d\tilde g(\gl).
\EA$$
Using $(\ref{B1})$ we obtain 
$$\BA {lll}
\myint{\{\abs{u_n}> t\}}{}\tilde g(\abs{u_n}) d\gs
\leq
c\norm\gm_{\frak M^b}^{\frac{\gth}{N-2}}\myint{t}{\infty}(\gl-M)^{-\frac{\gth}{N-2}} d\tilde g(\gl)
\\
\phantom{\myint{\{\abs{u_n}> t\}}{}\tilde g(\abs{u_n}) d\gs}
\leq\myfrac{c \gth}{N-2}\myint{t}{\infty}(\gl-M)^{-\frac{\gth}{N-2}-1}\tilde g(\gl) d\gl.
\EA$$
In view of assumption $(\ref{Z7})$, given $\ge>0$ we fix $t>M$ such that 
$$\myfrac{c \gth}{N-2}\myint{t}{\infty}(\gl-M)^{-\frac{\gth}{N-2}-1}\tilde g(\gl) d\gl
\leq \frac{\varepsilon}{2}.
$$
Then, setting $\gd=\frac{\ge}{2\tilde g(t)}$, we deduce 
$$\myint{E}{}d\gs\leq\gd\Longrightarrow \myint{E}{}\abs{g_2(u_n)}d\gs 
\le \varepsilon. 
$$
Since $g_1$ is bounded, this implies that $\{g(u_n)\}$ is uniformly integrable is $L^1_\gs(\Gw)$.
Since we already know that $g(u_n)\to g(u)$ $\sigma$-almost everywhere, 
it follows by  Vitali convergence's theorem that $g(u_n)\to g(u)$ in $L^1_\sigma(\Gw)$. 
Taking $\gz\in\BBX(\Gw)$ and letting $n\to\infty$ in the equality
$$-\myint{\Gw}{}u_n\Gd\gz dx+\myint{\Gw}{}g(u_n)\gz d\gs=\myint{\Gw}{}\gz d\gm_n
$$
yields the result.
\qeda

\mysection{The 2-D case}

In this section $\Gw$ is a bounded $C^2$ planar domain. The next result is the 2-D version of \rlemma{meas1}.
\blemma{meas2} Assume $N=2$ and $\gs\in \CM^+^+_{\frac{2}{2-\gth}}(\Gw)$ with $\gth>0$. If $\gm\in\frak M^b(\Gw)$ and $\gl\geq 0$,  we set
$E_\gl[\gm]:=\{x\in\Gw:\BBG[\abs\gm](x)>\gl\}$. Then 
\bel{2-D1}\BA {lll} \displaystyle
e_\gl^\gs(\gm):=\myint{E_\gl[\gm]}{}d\gs\leq \abs\Gw_\gs e^{1-\frac{\gl}{\gg\norm\gm_{\frak M^b}}}\qquad\text{for all }\gl>0,
\EA\ee
for some $\gg=\gg(\gth,{\rm diam}(\Gw))>0$
\es
\noindent\Proof If $\gm=\gd_a$ for some $a\in\Gw$, one has $0\leq\BBG[\gd_a](x)\leq \frac{1}{2\gp}\ln\left(\frac{d_\Gw}{\abs{x-a}}\right)$ where $d_\Gw={\rm diam}(\Gw)$. Hence 
$$E_\gl[\gd_a]\subset B_{d_\Gw e^{-2\gp\gl}}\Longrightarrow e_\gl^\gs(\gd_a)=\myint{E_\gl[\gd_a]}{}d\gs\leq cd^\gth_\Gw e^{-2\gth\gp\gl}.
$$
Let $E\subset\Gw$ be a Borel set, $\myint{E}{}d\gs=\abs E_\gs$ and $t>0$, then, as in \rlemma{meas1},
$$\BA {lll}\myint{E}{}\BBG[\gd_a]d\gs\leq t\myint{E}{}d\gs-\myint{t}{\infty}sde_s^\gs(\gd_a)\\[4mm]
\phantom{\myint{E}{}\BBG[\gd_a]d\gs}
\leq t\abs E_\gs+cd^\gth_\Gw\left(t+\myfrac{1}{2\gp\gth}\right)e^{-2\gth\gp t}.
\EA$$
If we choose $e^{-2\gth\gp t}=\frac{\abs E_\gs}{\abs \Gw_\gs}$ we infer
\bel{2-D2}\BA {lll} \displaystyle
\myint{E}{}\BBG[\gd_a]d\gs\leq \gg\abs E_\gs\left(\ln\left(\myfrac{\abs \Gw_\gs}{\abs E_\gs}\right)+1\right).
\EA\ee
For proving $(\ref{B1})$ we can assume that $\gm\geq 0$. Then there exists $\ga_j>0$ and $a_j\in\Gw$ such that
$$\gm=\sum_{j=1}^\infty\ga_j\gd_{a_j}\Longrightarrow\sum_{j=1}^\infty\ga_j=\norm\gm_{\frak M^b}.$$
Hence, for any Borel set $E\subset\Gw$, 
\bel{2-D3}\BA {lll}\displaystyle
\myint{E}{}\BBG[\gm](x)d\gs(x)=\sum_{j=1}^\infty\ga_j\myint{E}{}\BBG[\gd_{a_j}(x)]d\gs(x)
\leq \gg\abs E_\gs\left(\ln\left(\myfrac{\abs \Gw_\gs}{\abs E_\gs}\right)+1\right)\norm\gm_{\frak M^b}.
\EA\ee
If $E=E_\gl[\gm]$ we infer
$$\gl e_{\gl}^\gs(\gm)\leq \gg e_{\gl}^\gs(\gm)\left(\ln\left(\myfrac{\abs \Gw_\gs}{e_{\gl}^\gs(\gm)}\right)+1\right)\norm\gm_{\frak M^b},
$$
which implies the claim.\qeda

\bth{Vath1} Assume $N=2$, $\gs\in\CM^+^+_{\frac{2}{2-\gth}}(\Gw)$ with  $2\geq\gth >0$ and  $g:\BBR\mapsto\BBR$ a continuous function satisfying $(\ref{Z0})$. If  $a_\infty(g)=a_{\text{-}\infty}(g)=0$, for any $\gm\in\mathfrak M_b(\Gw)$ problem $(\ref{Z1})$ admits a very weak solution.
\es
\noindent\Proof Let $g^*$ be the monotone nondecreasing hull of $g$ defined by $(\ref{Z'9})$. If $m=\sup\{g(r):-r_0\leq r\leq r_0\}$ and $m'=\inf\{g(r):-r_0\leq r\leq r_0\}$ then $g\leq g^*+m$ on $\BBR_+$ and $g^*+m'\leq g$ on $\BBR_-$. If $\{\eta_n\}$ is a sequence of mollifiers and $\gm=\gm^+-\gm^-$, we set $\gm^+_n=\gm^+\ast\eta_n$, $\gm^-_n=\gm_-\ast\eta_n$, $\gm_n=\gm^+_n=-\gm^-_n$ and denote by $u_n$ the very weak solution of  
\bel{B8}\BA {lll}
-\Gd u_n+g(u_n)\gs=\gm_n\qquad&\text{in }\Gw\\
\phantom{-\Gd+g(u_n)\gs} 
u_n=0\qquad&\text{on }\prt\Gw.
\EA\ee
Since $\norm{\gm_n}_{L^1}\leq \norm{\gm}_{\mathfrak M_b}$, there holds by \rprop{Bre}, 
\bel{B9}\BA {lll}
\norm{u_n}_{L^1}+\norm{\gr g(u_n)}_{L_\gs^1}\leq c\norm{\gm}_{\mathfrak M_b}+M,
\EA\ee
and by \rlemma{L1},
\bel{B10}\BA {lll}
\norm{u_n}_{BMO}+\norm{\nabla u_n}_{L^{2,\infty}}\leq c\left(\norm{\gm}_{\mathfrak M_b}+\norm{\gr g(u_n)}_{L_\gs^1}\right)\leq c'\norm{\gm}_{\mathfrak M_b}.
\EA\ee
Again, there exists a set $E$ with $c_{1,q}(E)=0$ for any $q\leq 2-\gth$ such that $u_n(x)\to u(x)$ for all $x\in \Gw\setminus E$, hence 
$u_n(x)\to u(x)$ and $g(u_n(x))\to g(u(x))$ $d\gs$-almost everywhere in $\Gw$. This implies that $g(u)$ is $\gs$-measurable. In order to conclude we have to prove that $g(u_n)\to g(u)$ in $L^1_\gs(\Gw)$. Estimate $(\ref{2-D1})$ is valid, hence, for any $t>0$, 
$$\gt_n(t)=\myint{\{\abs{u_n(x)}>t\}}{}d\gs\leq e^\gs_{t-M}[\gm^+_n]+e^\gs_{t-M'}[\gm^-_n]\leq ce^{-\frac t{\gg \norm\gm_{\mathfrak M}}},
$$
by \rlemma{meas2}.  Since 
$$\abs{g(u_n)}\leq \left(g_+^*\left(u_n\right)-g_-^*\left(u_n\right)\right)+m-m',
$$
we have that
$$\BA {lll}\myint{E}{}\abs{g(u_n)}d\gs\leq \myint{E}{} g_+^*\left(u_n\right)d\gs-\myint{E}{} g_-^*\left(u_n\right)d\gs+(m-m')\abs E_\gs\\[4mm]
\phantom{\myint{E}{}\abs{g(u_n)}d\gs}
\leq -\myint{t}{\infty}g_+^*(s)d\abs{\{u_n>s\}}_\gs+\myint{-\infty}{-t}g_-^*(s)d\abs{\{u_n<s\}}_\gs+\left(m-m'\right)\abs E_\gs\\[4mm]
\phantom{\myint{E}{}\abs{g(u_n)}d\gs}
\leq -\myint{t}{\infty}\left(g_+^*(s)-g_-^*(-s)\right)d\gt_n(s)+\left(g_+^*(t)-g_-^*(-t)+m-m'\right)\abs E_\gs.
\EA$$
By integration by parts,
\bel{B11}\BA {lll}
-\myint{t}{\infty}\left(g_+^*(s)-g_-^*(-s)\right)d\gt_n(s)=\left(g_+^*(t)-g_-^*(-t)\right)\gt_n(t)+\myint{t}{\infty}\gt_n(s)d\left(g_+^*(s)-g_-^*(-s)\right)
\\[4mm]\phantom{-\myint{t}{\infty}\left(g_+^*(s)-g_-^*(-s)\right)d\gt_n(s)}
\leq \left(g_+^*(t)-g_-^*(-t)\right)\left(\gt_n(t)-ce^{-\frac{t}{\gg\norm\gm_{\frak M^b}}}\right)\\[4mm]
\phantom{------------------}+\myfrac{c}{\gg\norm\gm_{\frak M^b}}\myint{t}{\infty}e^{-\frac{s}{\gg\norm\gm_{\frak M^b}}}\left(g_+^*(s)-g_-^*(-s)\right)ds
\\[4mm]\phantom{-\myint{t}{\infty}\left(g_+^*(s)-g_-^*(-s)\right)d\gt_n(s)}
\leq \myfrac{c}{\gg\norm\gm_{\frak M^b}}\myint{t}{\infty}e^{-\frac{s}{\gg\norm\gm_{\frak M^b}}}\left(g_+^*(s)-g_-^*(-s)\right)ds.
\EA\ee
By assumption the integral on the right-hand side is convergent. We end the proof as in Theorem B, first by fixing $t$ large enough and then $\abs{E}_\gs$ small enough, and we derive the uniform integrability of  $\{g(u_n)\}$.\qeda\medskip

A similar result holds when $g$ has nonzero orders of growth at infinty. 
\bth{Vath2} Assume $N=2$, $\gs\in\CM^+^+_{\frac{2}{2-\gth}}(\Gw)$ with  $2\geq\gth >0$ and  $g:\BBR\mapsto\BBR$ a continuous function satisfying $(\ref{Z0})$. If  $0<a_\infty(g)<\infty$ and $-\infty<a_{\text{-}\infty}(g)<0$, there exists $\gd>0$ such that for any $\gm\in\mathfrak M_b(\Gw)$ satisfying 
$\norm\gm_{\mathfrak M_b}\leq \gd$ problem $(\ref{Z1})$ admits a very weak solution.
\es
\noindent\Proof The proof is a straightforward adaptation of the previous one. The choice of $\gd$ is such that 
\bel{B12}\BA {lll}
\norm\gm_{\mathfrak M_b}\leq \gd< \myfrac{1}{\gg}\sup\left\{\myfrac{1}{a_\infty(g)},-\myfrac{1}{a_{\text{-}\infty}(g)}\right\}
\EA\ee
and the conclusion follows from $(\ref{B11})$.\qeda


\mysection{The supercritical case}
\subsection{Proof of Theorem D}
{\it Proof of assertion I}. For $k>0$ set $g_k(r)=\max\{g(-k),\min\{g(k),g(r)\}\}$ and denote by $u_k$ the very weak solution of 
\bel{C0}\BA {lll}
-\Gd u+g_k(u)\gs=\gm&\qquad\text{in }\Gw\\
\phantom{-\Gd +g_k(u)\gs}
u=0&\qquad\text{on }\prt\Gw,
\EA\ee
which exists by Theorem B. 
 It follows from the proof of Theorem B (see $(\ref{Equ51})$ with $g=g_2$ and $g_1=0$) that 
\begin{equation}\label{C20} 
  \int_\Omega |g_k(u_k)|d\sigma \leq C, 
\end{equation} 
where the constant $C$ depends only on $\Omega$ and $|\mu|(\Omega)$. 
Thus the sequence of measures $\{g_k(u_k)\sigma\}$ is bounded.  
This implies that $\{u_k\}$ is bounded in 
$W^{1,q}(\Omega)$, $q<\frac{N}{N-1}$, and thus that, up to a subsequence, it converges in $L^1(\Omega)$ to some 
$u\in  W^{1,q}(\Omega)$, $q<\frac{N}{N-1}$. 
We can also assume that the convergence holds pointwise except on a set $E$ with zero $c_{1,q}$-capacity, which in turn is $\gs$-negligible by \rlemma{quasi} if we fix 
 $q\in \left(1,\frac{N}{N-1}\right)$ such that $\gth>N-q$. We also have that $u$ is finite but on a set with zero $c_{1,q}$-capacity hence $\gs$-negligible, therefore
 $$ g_k(u_k)\to g(u) \qquad \sigma\text{- \!almost everywhere.} $$ 
 Applying Fatou's lemma in 
$(\ref{C20})$ yields $g(u)\in L^1_\sigma(\Omega)$.  

By the maximum principle 
\bel{C1}\BA {lll}
-\BBG[\abs\gm]\leq u_k\leq \BBG[\abs\gm],
\EA\ee
hence
\bel{C2}\BA {lll}
g\left(-\BBG[\abs\gm]\right)\leq g_k(u_k)\leq g\left(\BBG[\abs\gm]\right),
\EA\ee
since $g$ is nondecreasing.

Because of assumption $(\ref{Z11})$ and in view of $(\ref{C2})$, we infer from Lebesgue dominated convergence that 
$\rho g_k(u_k)\to \rho g(u)$ in $L^1_\gs(\Gw)$. 
Thus we can pass to the limit in weak formulation of $(\ref{C0})$ with any $\zeta\in\BBX(\Gw)$.  

\medskip

\nind{\it Proof of assertion II}. We first notice that if $g$ is nondecreasing, vanishes at $0$ and satisfies $(\ref{Z12})$, then the function $g_k$ defined above also satisfies $(\ref{Z12})$ with the same constants $a$ and $b$.  We assume first that $\gm=\gm_r+\gm_s$ is nonnegative and we set $\gm^n_r=\gm_r\ast\eta_n$ where $\{\eta_n\}$ is a sequence of mollifiers. Let $u_k^n$ be the solution of $(\ref{C0})$ with right-hand side $\gm^n_r+\gm_s$ and $v_k^n$ the one of  $(\ref{C0})$ with right-hand side $\gm^n_r$ (in both cases existence and uniqueness
follows from Theorem B). Then $0\leq u_k^n\leq v_k^n+\BBG[\gm_s]$, $v_k^n\geq 0$ and 
$\BBG[\gm_s]\geq 0$. Since $g$ is non-decreasing, we deduce with $(\ref{Z12})$ that 
\bel{C4}\BA {lll}
0\leq g_k( u_k^n)\leq g_k\left(v_k^n+\BBG[\gm_s]\right)\leq a\left(g_k(v_k^n)+g_k(\BBG[\gm_s])\right)+b. 
\EA\ee
Since 
\bel{C5}\BA {lll}
\norm{v_k^n}_{L^1}+\norm{\gr g_k(v_k^n)}_{L_\gs^1}
\leq c\norm{\gm^n_r}_{\mathfrak M_b} 
\leq c\norm{\gm}_{\mathfrak M_b},
\EA\ee
up to subsequences, the sequences $\{v_k^n\}$ and $\{u_k^n\}$ converge in $L^1(\Gw)$ to some $v^n\in L^1(\Gw)$ and $u^n$ such that $\nabla v^n,\nabla u^n\in W^{1,q}$ for any $q<\frac{N}{N-1}$ when $k\to\infty$. As in I, $\{g_k(v_k^n)\}$ and $\{g_k(u_k^n)\}$ converge in $L_\gs^1(\Gw)$ to $\{g(v^n)\}$ and $\{g(u^n)\}$ respectively. Furthermore $v^n$ and $u^n$ satisfies 
\bel{C6}\BA {lll}
-\Gd v^n+g(v^n)\gs=\gm_r^n&\qquad\text{in }\Gw\\
\phantom{-\Gd +g(v^n)\gs}
v^n=0&\qquad\text{on }\prt\Gw,
\EA\ee
and 
\bel{C7}\BA {lll}
-\Gd u^n+g(u^n)\gs=\gm_s+\gm_r^n&\qquad\text{in }\Gw\\
\phantom{-\Gd +g(u^n)\gs}
u^n=0&\qquad\text{on }\prt\Gw,
\EA\ee
respectively and $0\leq u^n\leq v^n+\BBG[\gm_s]$. As in the proof of \rprop{Bre}, $ v^n\to v$ in $L^1(\Gw)$ and $\gr g(v^n)\to \gr g(v)$in $L_\gs^1(\Gw)$ as $n\to\infty$, and $v$ is a very weak solution of 
\bel{C8}\BA {lll}
-\Gd v+g(v)\gs=\gm_r&\qquad\text{in }\Gw\\
\phantom{-\Gd +g(v)\gs}
v=0&\qquad\text{on }\prt\Gw.
\EA\ee
As above $\{u^n\}$ converge in $L^1(\Gw)$ to some $u\in L^1(\Gw)$ (always up to some subsequence), there holds $u\leq v+\BBG[\gm_s]$ and $g( u^n)\to g(u)$ $\gs$-almost everywhere in $\Gw$ since the uniform bound on $\norm{\nabla u_n}_{L^{\frac{N}{N-1},\infty}}$ holds.  Furthermore
\bel{C9}\BA {lll}
0\leq g( u^n) \leq a\left(g(v^n)+g(\BBG[\gm_s])\right)+b\Longrightarrow
0\leq g( u) \leq a\left(g(v)+g(\BBG[\gm_s])\right)+b,
\EA\ee
and since $g(v^n)\to g(v)$ in $L^1_\gs(\Gw)$, the sequence $\{g( u^n)\}$ is uniformly integrable in  $L^1_\gs(\Gw)$. Again this implies that 
$g( u^n)\to g(u)$ in $L^1_\gs(\Gw)$ and $u$ is a very weak solution of $(\ref{Z1})$. \
If $\gm$ is signed measure, we construct successively the solutions $u^n_k$, $\overline u^n_k$ and $\underline u^n_k$  of $(\ref{C0})$ with right-hand side $\gm^n_r+\gm_s$, $\abs{\gm^n_r}+\abs{\gm_s}$ and $-\abs{\gm^n_r}-\abs{\gm_s}$ respectively, and the solutions $\overline v^n_k$ and $\underline v^n_k$ of $(\ref{C0})$ with right-hand side $\abs{\gm^n_r}$ and $-\abs{\gm^n_r}$ respectively.  Then
$\underline v^n_k-\BBG[\gm_s]\leq u_k^n\leq \overline v^n_k+\BBG[\gm_s]$ which implies by $(\ref{Z13})$
\bel{C10}\BA {lll}
a\left(g_k(\underline v^n_k)+g_k(-\BBG[\gm_s])\right)+b\leq g_k(u_k^n)\leq a\left(g_k(\overline v^n_k)+g_k(\BBG[\gm_s])\right)+b.
\EA\ee
Using the same estimates as above we conclude that $\displaystyle\lim_{n\to\infty}\lim_{k\to\infty}u_k^n=u$ exist in $L^1(\Gw)$, that 
$\displaystyle\lim_{n\to\infty}\lim_{k\to\infty}g_k(u_k^n)=g(u)$ holds $\gs$ almost everywhere in $\Gw$ and in $L^1_{\gs}(\Gw)$, which ends the proof.\qeda

\subsection{Reduced measures}
We adapt here some of the results in \cite{BMP} which turn out to be useful tools in our framework.
\blemma{monot} Let $\gs\in \CM^+^+_{\frac  N{N-\gth}}(\Gw)$ with $N\geq \gth> N-\frac{N}{N-1}$ and $g$ be nondecreasing satisfying $(\ref{Z0})$. Assume $\{\gm_n\}\subset\mathfrak M^+_b(\Gw)$ is an increasing sequence of good measures for problem $(\ref{Z1})$ converging to 
$\gm\in\mathfrak M^+_b(\Gw)$. Then $\gm$ is a good measure.
\es

\noindent\Proof Let $u_{\gm_n}$ be the solutions of $(\ref{Z1})$ with right-hand side $\gm_n$ then for any $n,k\in\BBN$, $ k\geq n$, we have 
since $u_0\in C^\ga(\overline\Gw)$, 
$$-m\leq u_0\leq u_{\gm_n}\leq u_{\gm_k} $$ 
for some $m\geq 0$ and then 
$$ g(-m)\leq g(u_0) \leq g(u_{\gm_n})\leq g(u_{\gm_k}). $$
We use $\zeta:=(\eta_1+\ge)^\ga-\ge^\ga$ as a test-function in the very weak formulation of the equation satisfied by $u_{\mu_n}-u_0$ as in the proof of $(\ref{AB61})$; then, recalling that $-\Delta\zeta\geq 0$, we obtain that 
$$ \int_\Omega (g(u_{\mu_n})-g(u_0))((\eta_1+\ge)^\ga-\ge^\ga]d\sigma
\leq \int_\Omega (\eta_1+\ge)^\ga d\mu_n 
\leq C\mu_n(\Omega)\leq C\mu(\Omega), $$ 
where $C$ is independent of $n$. 
letting successively $\ge\to 0$ and $\ga\to 0$ we obtain
$$ 0\leq \int_\Omega (g(u_{\mu_n})-g(u_0))d\sigma\leq C.$$ 
Hence $\{u_{\gm_n}\}$ is  bounded in $W^{1,q}_0(\Gw)$ for any $q<\frac{N}{N-1}$.
Thus there exists $u\in W^{1,q}_0(\Gw)$, $q<\frac{N}{N-1}$, such  
that $u_{\gm_n}\uparrow u$ in $L^1(\Omega)$ and pointwise but for a set $E$  with zero $c_{1,q}$-capacity. Since $\theta>N-\frac{N}{N-1}$ we can find some $q<\frac{N}{N-1}$ such that $\theta>N-q$. It then follows from 
\rlemma{quasi} that $\sigma(E)=0$.Thus $g(u_{\gm_n})\uparrow g(u)$ $\sigma$-almost everywhere.
Fatou's lemma yields $\int_\Omega (g(u)-g(u_0))d\sigma\leq C$, thus $g(u)\in L^1_\sigma(\Omega)$. By the dominated convergence theorem,  
$g(u_{\gm_n})\to g(u)$ in $L^1_\gs$. 
We can then pass to the limit in the equation satisfied by $u_{\gm_n}$ to obtain that $u=u_\mu$. 
\qeda

\bprop{VeronYarur} 
Assume $\gs$ and $g$ satisfy the assumptions of \rlemma{monot}.
Consider the set 
$$ Z = \left\{ x\in\Omega:\,  \int_\Omega \BBG(x,y)^q \rho(y)d\sigma(y) = \infty\right\}. $$ 
If $\mu\in\mathfrak{M}_b^+(\Omega)$ is such that $\mu(Z)=0$ then $\mu$ is good. 
\es

\nind\Proof 
We adapt to our case the proof of \cite{VeronYarur}[Thm 3.10].  
Consider the sets 
$$ C_n = \{ x\in\Omega:\,  \int_\Omega \BBG(x,y)^q \rho(y)d\sigma(y) \le n \}, 
\qquad n=1,2,\dots. $$ 
Since the function $x\to \int_\Omega \BBG(x,y)^q \rho(y)d\sigma(y) $ is lsc (by Fatou's lemma) the sets $C_n$ are closed. Moreover $C_n\subset C_{n+1}$ and 
$\bigcup_n C_n = \Omega\backslash Z$. 
Define $\mu_n:=1_{C_n}\mu$ i.e. $\mu_n$ is the measure $\mu$ restricted to $C_n$. 
Then each $\mu_n$ satisfies $(\ref{Z11})$. Indeed 
\begin{eqnarray*} 
\int_\Omega \BBG[|\mu_n|]^q\rho d\sigma 
& \le & \mu_n(\Omega)^{q-1}\int_\Omega\int_\Omega \BBG(x,y)^{q-1} d\mu_n(x)d\sigma(y)  \\
& \le & \mu(\Omega)^{q-1}\int_{C_n} \Big(\int_\Omega \BBG(x,y)^{q-1}d\sigma(y)\Big) 
 d\mu(x) \\ 
&\le & n\mu(\Omega)^q. 
 \end{eqnarray*}  
It follows from Theorem D that $\mu_n$ is good. 
Since $0\leq \mu_n\uparrow \mu$ we deduce from \rlemma{monot} that $\mu$ is good. 
\qeda

\blemma{abs} Assume $\gs$ and $g$ satisfy the assumptions of \rlemma{monot}. \smallskip

\nind I- If $\gm\in\mathfrak M^+_b(\Gw)$ is a good measure, any 
$\gn\in\mathfrak M^+_b(\Gw)$ such that $\gn\leq\gm$ is a good measure. \smallskip

\nind II- Let $\gm,\gm'\in\mathfrak M^+_b(\Gw)$. If $\gm$ and $-\gm'$ are good measures, any $\gn\in\mathfrak M_b(\Gw)$ such that 
$-\gm'\leq\gn\leq\gm$ is a good measure. 
\es

\noindent\Proof {\it Step 1}. Assume $\gm\in\mathfrak M^+_b(\Gw)$ is a good measure. For $k>0$ define $g_k$ by
$g_k(r)=\max\{g(-k),\min\{g(k),g(r)\}\}$, and denote by $u_{k,\gm}$ the solution of $(\ref{C0})$, which exists by Theorem B, and by $u_{\gm}$ the solutions of $(\ref{Z1})$.  Then $-m\leq u_0\leq \min\{u_\gm,u_{k,\mu}\}$. If $k>m$, then 
$g_k(u_{k,\gm})=\min\{g(k),g(u_{k,\gm})\}\leq g(u_{k,\gm})$. 
Hence
$$ -\Gd(u_\gm-u_{k,\mu})+\left(g_k(u_\gm)-g_k(u_{k,\gm})\right)\gs\leq 0.
$$
Then $u_\gm\leq u_{k,\mu}$ by \rlemma{Order}. Similarly $ u_{k',\mu}\leq u_{k,\mu}$ for $k'\geq k>m$. 
Using $\eta_1$ as test-function we obtain 
\bel{R1}\BA {lll}\myint{\Gw}{}(u_{k,\gm}-u_\gm)dx+\myint{\Gw}{}(g_k(u_{k,\gm})-g_k(u_\gm))\eta_1 d\gs=\myint{\Gw}{}(g(u_\gm)-g_k(u_\gm))\eta_1 d\gs.
\EA\ee
Since $g_k(r)\to g(r)$ for any $r\in\mathbb{R}$ and $|g_k(u_\gm)|\le |g(u_\gm)|$ with $\rho|g(u_\gm)|\in L^1_\sigma(\Gw)$, 
the right-hand side converges to $0$ as $k\to \infty$ and the second term on the left-hand side is nonnegative. Hence $u_{k,\gm}\to u_\gm$ in $L^1(\Gw)$ as $k\to \infty$, thus $\gr(g_k(u_{k,\gm})-g_k(u_\gm))\to 0$ in $L^1_\sigma(\Gw)$ which in turn yields 
$\gr g_k(u_{k,\gm})\to \gr g(u_\gm)$ in $L^1_\sigma(\Gw)$.
\smallskip

\nind {\it Step 2: proof of I}. Denote by $u_{k,\gn}$ the solution of 
\bel{R2}\BA {lll}
-\Gd u+g_k(u)=\gn&\qquad\text{in }\Gw\\
\phantom{-\Gd +g_k(u)}
u=0&\qquad\text{in }\prt\Gw.
\EA\ee
Then $-m\leq u_{k,\gn}\leq u_{k,\gm}$, $ u_{k',\mu}\leq u_{k,\mu}$ for $k'\geq k>m$ by \rlemma{Order}  and  $g_k(u_{k,\gn})\leq g_k(u_{k,\gm})$. 
Furthermore $\{u_{k,\gn}\}$ is bounded in $W^{1,q}_0(\Gw)$ for $1<q<\frac{N}{N-1}$ and
 thus relatively compact in $L^1(\Gw)$. 
Therefore there exists $u\in W^{1,q}_0(\Gw)$ such that $u_{k,\gn}\downarrow u$ 
in $L^1(\Gw)$ and also pointwise up to a set with zero $c_{1,q}$-capacity which is therefore a  $\gs$-negligible set. 
 By Step 1, the set 
$\gr g_k(u_{k,\gn})$ is uniformly integrable in $L_\gs^1(\Gw)$, this implies that $u=u_\gn$.
\smallskip

\nind {\it Step 3: Proof of II}. Because $-\gm'\leq \gn\leq\gm$ there holds $u_{k,-\gm'}\leq u_{k,\gn}\leq u_{k,\gm}$ and 
$g_k(u_{k,-\gm'})\leq g_k(u_{k,\gn})\leq g_k(u_{k,\gm})$. Since the sets $\{u_{k,-\gm'}\}$, $\{u_{k,\gn}\}$ and $\{u_{k,\gm}\}$ are relatively compact in $L^1(\Gw)$ and bounded in $W^{1,q}_0(\Gw)$ for $1<q<\frac{N}{N-1}$ and the sets $\{g_k(u_{k,-\gm'})\}$ and $\{g_k(u_{k,\gm})\}$ are uniformly integrable in $L_\gs^1(\Gw)$, then, up to a subsequence, $u_{k,\gn}\to u$ in $L^1(\Gw)$ and $\gs$-almost everywhere as $k\to\infty$. This implies that  $g(u)\in L^1_\gs(\Gw)$ and 
$\gr g_k(u_{k,\gn})\to \gr g(u)$ in $L^1_\gs(\Gw)$. Hence $u=u_\gn$. \qeda

\medskip

The proof of the next result, based upon Zorn's lemma, is a variant of the one of \cite[Th 4.1]{BMP} which uses inverse maximum principle 
\cite[Corollary 4.8]{BMP}.

\blemma{redu} Assume $\gs$ and $g$ satisfy the assumptions of \rlemma{monot}. If $\gm\in\mathfrak M^+_b(\Gw)$ there exists a largest 
good measure smaller than $\gm$, and it is nonnegative. 
\es

\noindent\Proof 
Let $\CZ_\gm$ be the subset of all bounded nonnegative good measures smaller than $\gm$. 
Notice first that $\CZ_\gm$ is non-empty since it contains the regular part $\gm_r$ of $\gm$  with respect to the N-dim Hausdorff measure.  
We now show that $\CZ_\gm$ is inductive. Let $\CC_I:=\{\gm_i\}_{i\in I}$ be a totally ordered subset of $\CZ_\gm$. For $\gz\in C_0(\overline\Gw)$, $\gz\geq 0$, the set of nonnegative real numbers
$$\CC_I(\gz):=\left\{\myint{\Gw}{}\gz d\gm_i\right\}
$$
is bounded from above by $\myint{\Gw}{}\gz d\gm$. Note that can we extend $\gm$ as a positive linear form on $C_0(\overline\Gw)$ since it is a Radon measure and $\gm(\prt\Gw)=0$.  Hence $\CC_I(\gz)$ admits an upper bound $L(\gz)$ and there exists a sequence $\{i_k\}\subset I$ such that 
$$\myint{\Gw}{}\gz d\gm_{i_k}\uparrow L(\gz)\leq \myint{\Gw}{}\gz d\gm\qquad\text{as } k\to\infty. 
$$
By the Stone-Weiertrass theorem there exists a dense subset $\{\gz_n\}$ of the set of nonnegative elements in $C_0(\overline\Gw)$. By Cantor diagonal process there exists a subsequence $\{i_{n_k}\}\subset I$ such that 
$$\myint{\Gw}{}\gz_n d\gm_{i_{n_k}}\uparrow L(\gz_n)\leq \myint{\Gw}{}\gz_n d\gm\qquad\text{as } k\to\infty.
$$
Clearly the map $\gz_n\mapsto L(\gz_n)$ is additive, positively homogeneous of order one and  satisfies
$$L(\gz)\leq \myint{\Gw}{}\gz d\gm\qquad\text{for all }\gz\in C_0(\overline\Gw),\,\gz\geq 0.
$$
Hence $L$ extends as a positive linear functional on $C_0(\overline\Gw)$, dominated by $\gm$ denoted by $\gm_{\CC_I}$. Since $\gm$ is a Radon measure in $\Gw$, $\gm_{\CC_I}(\prt\Gw)=0$, hence it is a Radon mesure. Furthermore it is a good measure by \rlemma{monot}. 
It follows that $\gm_{\CC_I}\in \CZ_\gm$ . 
Moreover since $L(\zeta)$ is an upper bound of $\CC_I(\gz)$ for any nonegative 
$\zeta\in C_0(\overline\Gw)$, we have $\gm_{\CC_I}\geq \mu_i$ for any $i\in I$. 
Hence the set $\CZ_\gm$ is inductive. 

 As a consequence of Zorn's lemma,  $\CZ_\gm$ admits at least one maximal element that we denote $\gm^*$. If $\gn$ is any nonnegative good measure smaller than $\gm$ it belongs to  $\CZ_\gm$ and hence it cannot dominate $\gm^*$. 
 It remains to prove that $\gn\leq \gm^*$. 
 Set $\gl=\sup\{\gn,\gm^*\}$ and let $\gl^*$ be a maximal element of $\CZ_{\gl}$. 
 Since $\nu$ and $\mu^*$ are good measures, we have $\gn^*=\gn$ and 
 $(\gm^*)^*=\gm^*$. It follows that  
$\gl^*\geq\gn^*=\gn$ and $\gl^*\geq (\gm^*)^*=\gm^*$ so that  
$\gl^*\geq \sup\{\gn,\gm^*\}=\gl$. This implies that $\gl^*=\gl\geq\gm^*$. 
On the other hand, since $\nu,\mu^*\le \mu$, we have $\lambda\le \mu$ and thus 
 $\lambda^*\le \mu$. By definition of a maximal element it implies  that $\gl^*=\gl=\gm^*$, and finally $\gm^*=\sup\{\gn,\gm^*\}$. 
We infer $\gn\leq \gm^*$ and then $\gm^*$ is the maximum of $\CZ_\gm$.\qeda

\bcor{sup} Assume $\gs$ and $g$ satisfy the assumptions of \rlemma{monot}. If $\gm,\gn \in\mathfrak M^+_b(\Gw)$ are good measures, then 
$\sup\{\gm,\gn\}$ is a good measure.
\es
\noindent\Proof Set $\gl=\sup\{\gm,\gn\}$. Then 
\bel{R3}\BA {lll}\gl\geq \gl^*=(\sup\{\gm,\gn\})^*\geq \sup\{\gm^*,\gn^*\}=\sup\{\gm,\gn\}=\gl.
\EA\ee
This implies $\gl= \gl^*$, hence $\gl$ is a good measure. \qeda

\subsection{The capacitary framework}
We start with the following regularity estimate for the Poisson problem

\blemma{regul} For any $s\geq 0$ and $1<p<\infty$, the mapping $\gm\mapsto \BBG[\gm]$ is continuous from $\mathfrak M_b(\Gw)\cap H^{s-2,p}(\Gw)$ to  $H^{s,p}(\Gw)$.
\es
\Proof It is classical that the mapping $G_D:\gl\mapsto u=G_D(\gl)$ solution of $-\Gd u=\gl$ in $\Gw$ and $u=0$ on $\prt\Gw$ is continuous from  $H^{s-2,p}(\Gw)$ to  $H^{s,p}(\Gw)$ for $1<p<\infty$ and $s>\frac 1p$ (see e.g. \cite[Example 3.15 p. 314]{G}). Thus we are left with the case 
$0\leq s\leq \frac 1p$. If $\gl\in \mathfrak M_b(\Gw)$, then $G_D(\gl)=\BBG[\gl]$ is a very weak solution, hence, since $\BBX(\Gw)\subset \displaystyle C^1_c(\overline\Gw)\cap\left(\bigcap_{1<r<\infty}H^{2,r}(\Gw)\right)$,
$$-\myint{\Gw}{}G_D(\gl)\Gd\gz dx=\myint{\Gw}{}\gz d\gl\leq \norm{\gz}_{H^{2-s,p'}}\norm{\gl}_{H^{s-2,p}}\quad\text{for all }\,\gz\in \BBX(\Gw).
$$
In particular, if $\gz=\BBG[v]$, then  $\norm{\gz}_{H^{2-s,p'}}\leq c\norm{v}_{H^{-s,p'}}$ since $-s>-2+ 1/p'$, and 
$$\myint{\Gw}{}G_D(\gl)v dx\leq c\norm{v}_{H^{-s,p'}}\norm{\gl}_{H^{s-2,p}}\quad\text{for all }\,v\in \Gd(\BBX(\Gw)).
$$
In particular this inequality holds if $v\in C_c(\overline\Gw)$ which is dense in $H^{-s,p'}(\Gw)$. Finally this inequality means that the mapping 
$v\mapsto\myint{\Gw}{}G_D(\gl)v dx$ is a continuous linear form over $H^{-s,p'}(\Gw)$, it thus belongs to $H^{s,p}(\Gw)$.\qeda

\bprop{E} Let $\gs$ and $g$ satisfy the assumptions in Theorem E. If $\gm\in\mathfrak M_b(\Gw)$ is such that $\abs\gm\in H^{s-2,p}(\Gw)$ for some $p>1$ and $s>0$ such that 
$N-\gth<sp<N$ and $\frac{\gth p}{N-sp}\geq q$, then $(\ref{Z2})$ admits a unique very weak solution.
\es
\noindent\Proof By \rlemma{regul}, if $\abs\gm\in H^{s-2,p}(\Gw)$ then $\BBG[\abs\gm]\in H^{s,p}(\Gw)$. 
By \rprop{AdamsMazya} 
$$\norm{\BBG[\abs\gm]}_{L^q_\gs}\leq c\norm{\BBG[\abs\gm]}_{H^{s,p}}
$$
if and only if $\gs\in \CM^+^+_{r}(\Gw)$ with $\frac 1r= q\left(\frac 1q-\frac 1p+\frac sN\right)=\frac{N-\gth'}{N}$. Then $q=\frac{\gth'p}{N-sp}$. Hence, if 
$\frac{\gth p}{N-sp}\geq q$ we get $\gth\geq \gth'$ and then $\CM^+^+_{\frac{N}{N-\gth}}(\Gw)\subset \CM^+^+_{\frac{N}{N-\gth'}}(\Gw)$ 
by $[\ref{Y7})$. We conclude by Theorem D.\qeda\medskip

\nind{\Remark} This result 
covers the case $q=p$, in which  any bounded measure such that $\abs\gm\in  H^{\frac{N-\gth}{q}-2,q}(\BBR^N)$ is eligible for solving problem $(\ref{Z1})$. \medskip

\medskip

\nind{\it Proof of Theorem E}. If $\gm$ is absolutely continuous with respect to the $c_{2-s,p'}$-capacity, so are $\gm^+$ and $-\gm^-$. By \cite{FePr} there exists an increasing sequence of positive bounded Radon measures $\gm_j\in H^{s-2,p}(\Gw)$ converging to $\gm^+$. By \rprop{E} $\gm_j$ is a good measure, hence by \rlemma{monot} $\gm^+$ is a good measure. In the same way $-\gm^-$ is a good measure. Since 
$-\gm_-\leq \gm\leq \gm_+$, it follows from \rlemma{abs}-II that $\gm$ is a good measure.\qeda\medskip

\medskip

\nind{\it Proof of \rprop{CondSimple}}. 
 Notice first that if $\mu \in \CM^+_{\frac{N}{N-\gth^*}}(\Gw)$ with 
 $\theta^*>N-sp$, then for any compact $K\subset \Omega$, 
 \begin{equation}\label{CondSimple1} 
  |\mu|(K)\leq c' \left(c_{(s,p)}(K)\right)^\frac{1}{p}. 
  \end{equation} 
  In particular $\mu$ is absolutely continuous w.r.t $c_{(s,p)}$-capacity. 
  Indeed under the assumption on $\gth^*$ we have 
  $H^{s,p}(\Omega)\hookrightarrow L^1_{|\mu|}(\Omega)$. 
  It follows that for any $v\in H^{s,p}(\Omega)$, $v\geq 1$ on $K$, we have 
  $$ |\mu|(K)\leq \int_Kv d|\mu| \leq \|v\|_{L^1_{|\mu|}}
  \leq C\|v\|_{H^{s,p}}. $$
  We deduce $(\ref{CondSimple1})$ taking the infimum over $v$. 
  To apply Theorem E we need $\mu$ to be $c_{2-\frac{N-\theta}{q},q'}$-diffuse. It thus suffices to take $\gth^*>N-sp$ with 
  $s=2-\frac{N-\theta}{q}$ and $p=q'$. 
  We obtain exactly the condition on $\theta^*$ stated in 
  \rprop{CondSimple}. \qeda

\subsection{The case $g(u)=\abs u^{q-1}u$.}

In the sequel we consider the following equation
    \bel{D1}\BA {lll} 
-\Gd u+\abs{u}^{q-1}u\gs=\gm & \qquad\text{in }\Gw\\
\phantom{-\Gd +\abs{u}^{q-1}u\gs}
u=0 & \qquad\text{in }\prt\Gw,
\EA\ee
where $q>1$. A measure for which there exists a solution, necessarily unique  by \rlemma{Uniq}, is called {\it q-good}. Assume that  $\gs\in \CM^+^+_{\frac{N}{N-\theta}}$ with 
$N\geq \theta>N-\frac{N}{N-1}$.  Then the critical exponent $q$ from the point of view of
 $(\ref{Z7})$ in Theorem B is
\bel{D2^*}\BA {lll} 
q_\gth:=\myfrac{\gth}{N-2},
\EA\ee
which is larger than $1$ if $N>2$. \smallskip
 
Let $q>1$ and  $\sigma\in\mathfrak M_b^+(\Gw)$. Recall that the Green function $G$ of the Dirichlet Laplacian in $\Gw$ is defined on 
$\overline\Gw\times \overline\Gw$ with values in $[0,+\infty]$ with $G(x,x)=+\infty$, $x\in\Gw$, and $G(x,y)=0$ if $x\in \partial\Gw$ or  $y\in \partial\Gw$.
We extend $G$ to $\BBR^N\ti\overline\Gw$ by setting $G(x,y)=0$ if $(x,y)\in \overline\Gw^c\ti\overline\Gw$. Hence $x\mapsto G(x,y)$ is lower semicontinuous in $ \BBR^N$ and 
$y\mapsto G(x,y)$ is lower semicontinuous in $\Gw$,  and thus is $\gs$-measurable. 
Following \cite[Sec. 2.3]{AdHe} we then consider the following set function with value in $[0,+\infty]$, 
\bel{D20'}\BA {lll} 
c_{q}^\gs(E)
=\inf\left\{\myint{\Gw}{}\abs v^{q'}d\gs: v\in L_\gs^{q'}(\Gw),\,\BBG[v\gs](x)\geq 1 \;\,\text{for all } x\in E\right\},
\EA\ee
for any $E\subset \Gw$. 
According to the general theory developped in \cite[Sec. 2.3]{AdHe} $c_{q}^\gs$ is a regular capacity in the sense of Choquet. Using  the lower semicontinuity of $y\mapsto \BBG[v\gs](y)$ 
(see\cite[Prop 2.3.2]{AdHe}) it is easy to verify that for any compact set $K\subset\Gw$, there holds
\bel{D3-0}\BA {lll} 
c_{q}^\gs(K)
=\inf\left\{\myint{\Gw}{}\abs v^{q'}d\gs: v\in L_\gs^{\infty}(\Gw),\,\BBG[v\gs](x)\geq 1 \;\,\text{for all } x\in K\right\}.
\EA\ee

The dual formulation of the capacity is the following (see \cite[Th 2.5.1]{AdHe}),
   \bel{D3}\BA {lll} 
\left(c_{q}^\gs(K)\right)^{\frac{1}{q'}}=\sup\left\{\gl(K): \gl\in\mathfrak M_b^+(K),\, 
\norm {\BBG[\gl]}_{L^q_\gs}\leq 1\right\}\quad \text{for $K\subset\Gw, K$ compact}.
\EA\ee
Existence of extremal measures satisfying equality in $(\ref{D3})$ is proved in \cite[Th 2.5.3]{AdHe}. \medskip

\nind\Remark Note that the $\geq$ inequality in $(\ref{D3})$ follows directly from the following one 
\bel{D31}\BA {lll} 
\nu(K)\le \left(c_q^\sigma(K)\right)^{\frac{1}{q'}} \|\BBG\nu\|_{L^q_\sigma},
\EA\ee
which holds for any $\gn\in\mathfrak M_b^+(\Omega)$ such that $\BBG[\nu]\in L^q_\sigma$ 
and any $K\subset\Omega$ compact. 
\medskip

We now give some sufficient conditions for a bounded measure to be absolutely continuous with respect to the capacity $c_{q}^\gs$. First in view of $(\ref{D31})$ and the dual expression of the capacity it is clear that there holds:

\blemma{capa0}  
If $\nu\in\mathfrak M_b(\Gw)$ is such that $\BBG[\abs\nu]\in L^q_\gs(\Gw)$, then 
$\nu$ is absolutely continuous with respect to the capacity $c_{q}^\gs$.
This holds in particular if 
$\gn\in\mathfrak M_b(\Gw)$ is such that 
$\abs\gn\in H^{s-2,p}(\Gw)$ for some $p>1$ and $s>0$ verifying 
$N-\gth<sp<N$ and $\frac{\gth p}{N-sp}\geq q$.  
\es

As a direct consequence we have 

\blemma{capa3}  
If $\nu\in\mathfrak M_b(\Gw)$ is $c_{2-s,p'}$-diffuse where� $s$ and $p$ are as in 
\rlemma{capa0}, then $\nu$ is absolutely continuous with respect to the capacity $c_q^\gs$.
\es

\noindent\Proof 
If $\nu\geq 0$ there exists a sequence of nonnegative measures 
$\{\nu_n\}\subset H^{s-2,p}(\Gw)$ such that $\nu_n\uparrow \nu$. 
If $K$ is a compact such that $c_q^\sigma(K)=0$ then $\nu_n(K)=0$ by 
\rlemma{capa0} and thus $\nu(K)=0$. When $\nu$ is a signed measure, we apply the above to $\nu^\pm$. 
\qeda

The following particular case will be useful: 

\blemma{capa1}  
If $\nu\in\CM^+_{\frac{N}{N-\gth}}(\Gw)$ with $N\geq \gth>N-2$,  then 
$\nu$ is absolutely continuous with respect to the capacity $c_{q}^\gs$. \es

\noindent\Proof 
We have $|\nu|\in \CM^+_p(\Gw)$ for some $p>\frac{N}{2}$. We then obtain from 
$(\ref{Y2+3})$ that $\BBG[|\gn|]$ is bounded so that 
$\BBG[|\gn|]\in L^q_\gs(\Gw)$. The conclusion follows from the previous lemma. 
\qeda

\medskip

\nind\Remark It is noticeable that if the support of a nonnegative measure $\gm$ does not intersect the support of $\gs$, it is always $q$-good. This is due to the fact that $\BBG[\gm]$ is bounded on the support of $\gs$, hence $\BBG[\gm]\in L^q_\gs(\Gw)$ for any $q<\infty$ and the result follows from Theorem D. Hence, a more accurate necessary condition must involve a notion of density of $\gs$ on its support,  a property which has been developed by Triebel \cite{Tri4} in connection with fractal measures. \medskip

We recall that the $\gth$-dimensional Hausdorff measure $H^\gth$, $0\leq \gth\leq N$, is defined on subsets $E$ of $\BBR^N$ by
\bel{H1}\BA {lll} \displaystyle
 H^\gth (E)=\lim_{\gd\to 0}\left(\inf\left\{\sum_{j=1}^\infty(\text{diam}\,U_j)^\gth: E\subset\bigcup_{j=1}^\infty U_j,\text{diam}\,U_j\leq\gd\right\}\right).
 \EA\ee

\bdef{regular} A nonnegative Radon measure $\gs$ on $\overline\Gw$ with support $\Gg$ is $\gth$-regular with $0\leq \gth\leq N$ if there exists $c>0$ such that 
  \bel{H2}\BA {lll} \displaystyle
\frac{1}{c} r^\gth\leq\abs{B_r(x)}_\gs\leq cr^\gth\qquad\qquad\text{for all } x\in\Gg\,,\;\text{for all } r>0.
 \EA\ee
 The support $\Gg$ of $\gs$ is called a $\gth$-set.
  \es
  
By \cite[Th 3.4]{Tri4} $\gs$ is equivalent in $\overline\Gw$ to the restriction  $H^\gth\lfloor_\Gg$ of $H^\gth$ to $\Gg$ in the sense  that there exists $c'>0$ such that 
  \bel{H3}\BA {lll} \displaystyle
\frac{1}{c'} H^\gth (E\cap\Gg)\leq \gs(E)\leq c'H^\gth (E\cap\Gg)\qquad\text{for all } E\subset\overline\Gw\,, \;E\text{ Borel}.
 \EA\ee

The description of $L^p_\gs(\Gg)$ necessitates to introduce the scale of Besov spaces and their {\it trace} on $\Gg$. For $0<s<1$, $1\leq p,q\leq\infty$, we denote by $B^s_{p,q}(\Gw)$ the space obtained by the real interpolation method by
  \bel{H4}\BA {lll} \displaystyle
B^{s}_{p,q}(\Gw)=\left [W^{1,p}(\Gw),L^{p}(\Gw)\right ]_{s,q}.
 \EA\ee
 Details can be found in \cite{Tri1}. It's norm is equivalent to 
  \bel{H5}\BA {lll} \displaystyle
\norm\gf_{B^{s}_{p,q}}=\norm  v_{L^p}+\left(\myint{0}{\infty}\myfrac{\left(\gw_p(t;  v)\right)^q}{t^{sq}}\myfrac{dt}{t}\right)^{\frac 1q},
 \EA\ee
 if $q<\infty$ and 
   \bel{H6}\BA {lll} \displaystyle
\norm\gf_{B^{s}_{p,\infty}}=\norm  v_{L^p}+\sup_{t>0}\myfrac{\gw_p(t;  v)}{t^{s}}, \EA\ee
 where
 $$\gw_p(t;\gf)=\sup_{\abs h<t}\norm{  v(.+h)-  v(.)}_{L^p}
 $$
 For $k\in\BBN_*$, $B^{k+s}_{p,q}(\Gw)=\{  v\in W^{k,p}(\Gw):D^\ga  v\in B^{s}_{p,q}(\Gw)\,,\;\text{for all }\ga\in\BBN^N,\,\abs\ga=k\}$ with norm 
 $$\norm  v_{B^{k+s}_{p,q}}=\norm  v_{W^{k-1,p}}+\sum_{\abs\ga=k}\norm{D^\ga  v}_{B^{s}_{p,q}}.
 $$
 
 If $\Gg\subset\BBR^N$ is a closed set with zero Lebesgue measure, 
    \bel{H8}\BA {lll} \displaystyle
 B^{s,\Gg}_{p,q}(\BBR^N)=\left\{  v\in  B^{s}_{p,q}(\BBR^N):\langle v,\gf\rangle=0\quad\text{for all }\gf\in\CS(\BBR^N)\text{ s.t. }\gf\lfloor_\Gg =0\right\},
  \EA\ee
  where
  $$\langle v,\gf\rangle=\myint{\BBR^N}{}v\gf dx,
  $$
  is the pairing between $\CS'(\BBR^N)$ and $\CS(\BBR^N)$.
 If  $v\in L^q_\gs(\Gw)$ and $\gs$ has  support $\Gg\subset\overline\Gw$, the linear map
   \bel{H9}\BA {lll} \displaystyle
\gf\mapsto T^\gs_v(\gf)=\myint{\Gg}{}\gf vd\gs
 \EA\ee
defined on $\CS(\BBR^N)$ is a tempered distribution in $\BBR^N$. The following results are proved in \cite[Th 18.2, 18.6]{Tri4}.
\bprop {Tr1}Assume $\gs$ is $\gth$-regular, $0<\gth<N$ with support $\Gg\subset\BBR^N$. Then for any $1< p\leq\infty$ the mapping 
$v\mapsto T_v^\sigma$ satisfies
   \bel{H10}\BA {lll} \displaystyle
\abs{T^\gs_v(\gf)}\leq c\norm v_{L^p_\gs}\norm \phi_{B^{\frac{N-\gth}{p'},\Gg}_{p',1}}\qquad\text{for all } \gf\in \CS(\BBR^N).
 \EA\ee
 Furthermore this mapping is onto, that we write $L^p_\gs(\Gg)\sim \left(B^{\frac{N-\gth}{p'},\Gg}_{p',1}\right)'=B^{-\frac{N-\gth}{p'},\Gg}_{p,\infty}$.
\es

\bprop {Tr2} Assume $\gs$ is $\gth$-regular, $0<\gth<N$ with support $\Gg\subset\BBR^N$. Then for any $1< p\leq\infty$ the restriction operation
from $\CS(\BBR^N)$ to $C(\Gg)$, $\phi\mapsto \phi\lfloor_\Gg$ can be extended as a continuous linear operator from  $B^{\frac{N-\gth}{p}}_{p,1}(\BBR^N)$ to $L^p_\gs(\Gg)$ that we denote $Tr_\Gg$. Furthermore this operator is onto.
\es

\bdef{captri} If $\gs\in \mathfrak M_b^+(\Gw)$ is $\gth$-regular, $N\geq \gth>N-2$ with support $\Gg\subset\Gw$ and $m,q>1$, we set
   \bel{H11}\BA {lll} \displaystyle
c^{2-\frac{N-\gth}{q},\Gg}_{q',\infty}K)=\inf\left\{\norm\gz^{q'}_{B^{2-\frac{N-\gth}{q}}_{q',\infty}}:\gz\in B^{2-\frac{N-\gth}{q},\Gg}_{q',\infty}(\Gw)\text{ s.t. }\,\gz\geq \chi_{_K} \right\},
 \EA\ee
 where
   \bel{H12}\BA {lll}B^{2-\frac{N-\gth}{q},\Gg}_{q',\infty}(\Gw)=\left\{\gz\in B^{2-\frac{N-\gth}{q}}_{q',\infty}(\Gw)\text{ s.t. }\Gd\gz\in B^{-\frac{N-\gth}{q},\Gg}_{q',\infty}(\Gw)\right\}.
    \EA\ee
    Notice that $B^{2-\frac{N-\gth}{q},\Gg}_{q',\infty}(\Gw)$ is a closed subspace of $B^{2-\frac{N-\gth}{q}}_{q',\infty}(\Gw)$.
\es

\bprop {comparison} Assume $\gs\in \mathfrak M_b^+(\Gw)$ is $\gth$-regular, $N\geq \gth>N-2$ with support $\Gg\subset\Gw$ and $q>1$. Then there exists a positive constant
$M>0$ such that 
    \bel{H13}\BA {lll} \displaystyle
\myfrac 1Mc^\gs_q(K)\leq c^{2-\frac{N-\gth}{q},\Gg}_{q',\infty}(K)\leq Mc^\gs_q(K),
 \EA\ee
for all compact set $K\subset\Gw$.
\es
\noindent\Proof By standard elliptic equations and interpolation theory (see \cite{Tri1}, \cite{Tri2}), for any $\psi\in B^{-\frac{N-\gth}{q},\Gg}_{q',\infty}(\Gw)$, $\BBG[\psi\gs]\in B^{2-\frac{N-\gth}{q}}_{q',\infty}(\Gw)$ and there holds
    \bel{H14}\BA {lll} \displaystyle
\myfrac 1c\norm{\BBG[\psi\gs]}_{B^{2-\frac{N-\gth}{q}}_{q',\infty}}\leq \norm{\psi}_{B^{-\frac{N-\gth}{q},\Gg}_{q',\infty}}\leq 
c\norm{\BBG[\psi\gs]}_{B^{2-\frac{N-\gth}{q}}_{q',\infty}}.
 \EA\ee
 By \rprop{Tr1} we can replace $\norm{\psi}_{B^{-\frac{N-\gth}{q},\Gg}_{q',\infty}}$ by $\norm{\psi}_{L^{q'}_\gs}$ in the above inequality, up to a change of constants $c$. Let $\{v_k\}\subset L_\gs^{\infty}(\Gw)$ such that $v_k\geq 0$, $\gz_k:=\BBG[v_k\gs]\geq 0$ on $K$ and $\norm{v_k}_{L^{q'}_\gs}\downarrow \left(c^\gs_q(K)\right)^{\frac 1{q'}}$. Since $(\ref{H12})$ is equivalent to
 $$ \myfrac 1c\norm{\gz_k]}_{B^{2-\frac{N-\gth}{q}}_{q',\infty}}\leq \norm{v_k}_{L^{q'}_\gs}\leq 
c\norm{\gz_k]}_{B^{2-\frac{N-\gth}{q}}_{q',\infty}},
 $$
 we derive $c^{2-\frac{N-\gth}{q},\Gg}_{q',\infty}(K)\geq \frac{1}{c^{q'}}c^\gs_q(K)$. Similarly $c^{2-\frac{N-\gth}{q},\Gg}_{q',\infty}(K)\leq c^{q'}c^\gs_q(K)$.\qeda\medskip
 
 \nind{\it Proof of Theorem F}. By \rlemma {capa1} the measure $u^q$ vanishes on Borel sets with zero $c^\gs_q$-capacity. Since $u\in L^q_\gs(\Gw)$ the mapping 
 $$\gf\mapsto=\myint{\Gg}{}u\gf d\gs=\langle u,\gf\rangle
 $$
 is a tempered distribution that we denote by $T_u^\gs$, hence
 $$\abs{\langle \Gd u,\gf\rangle}=\abs{\langle  u,\Gd\gf\rangle}=\abs{\myint{\Gw}{}u\Gd\gf d\gs}\leq 
 \norm u_{L^q_\gs} \norm {\Gd\gf}_{L^{q'}_\gs}.
 $$
 Using \rprop{Tr1} 
 $$\norm {\Gd\gf}_{L^{q'}_\gs}\leq c\norm {\Gd\gf}_{B^{-\frac{N-\gth}{q},\Gg}_{q',\infty}}\leq c'\norm {\gf}_{B^{2-\frac{N-\gth}{q},\Gg}_{q',\infty}}.
 $$
 Therefore the nonnegative measure $T_u^\gs$ is a continuous linear form on $B^{2-\frac{N-\gth}{q},\Gg}_{q',\infty}(\Gw)$. Therefore it vanishes on 
 Borel sets with zero $c^{2-\frac{N-\gth}{q},\Gg}_{q',\infty}$-capacity, which actually coincide with Borel sets with zero zero $c^\gs_q$-capacity.\qeda

\subsection{Removable singularities}

It is easy to prove that for any compact set $K\subset\Gw$, there exists $\gm_K\in\mathfrak M_b^+(K)$ such that $\myint{\Gw}{}(\BBG[\gm_K])^qd\gs= 1$ and 
$c_q^\gs(K)=\gm_K(K)$ (see \cite{AdHe}[Th 2.5.3]). 
Since $\gm_K$ is an admissible measure, it follows from Theorem D that  $(\ref{Z2})$ is solvable with $\gm=\gm_K$, hence $K$ is not removable. Although it could be conjectured that a compact set with zero $c_q^\gs$-capacity is removable we can prove this assertion only for sigma-moderate solutions.

\bdef{sigma} Let $q>1$, $\gs\in\CM^+^+_{\frac{N}{N-\gth}}(\Gw)$ where $N\geq\gth>N-2$ and $K\subset \Gw$ a compact set. A nonnegative function $u\in L^1_{loc}(\overline\Gw\setminus K)\cap L^q_{\gs,\,loc}(\overline\Gw\setminus K)$ 
is a sigma-moderate solution of
    \bel{E1}\BA {lll} 
-\Gd u+\abs{u}^{q-1}u\gs=0\qquad&\text{in }\Gw\setminus K\\
\phantom{-\Gd +\abs{u}^{q-1}u\gs}
u=0\qquad&\text{in }\prt\Gw,
\EA\ee
if there exists an increasing sequence $\{\gm_n\}\subset\mathfrak M_b^+(K)$ of $q$-good measures such that $u_{\gm_n}\to u$ in 
$L^1_{loc}(\overline\Gw\setminus K)\cap L^q_{\gs\,loc}(\overline\Gw\setminus K)$.
\es \smallskip

\bth{remov1} Under the assumptions on $q$, $\gs$ and $K$ of \rdef{sigma}, if $c_q^\gs(K)=0$ then the only sigma-moderate solution of 
$(\ref{E1})$ is trivial.
\es
\noindent\Proof Since $c_q^\gs(K)=0$ the set of nonnegative $q$-good measures with support in $K$ is reduced to the zero function by Theorem F. This implies the claim. \qeda\medskip

\nind\Remark We conjecture that for any compact set $K\subset\Gw$, any nonnegative local solution of $(\ref{R1})$ is sigma-moderate. This would imply that a necessary and sufficient condition for a local nonnegative solution of $(\ref{R1})$ to be a solution in $\Gw$ is $c_q^\gs(K)=0$. However this type of result is usually difficult to prove, see \cite{Mse}, \cite{MaVe2}, \cite{Dyn} in the framework of semilinear equations with measure boundary data. \medskip

In order to find necessary and sufficient conditions for the removability of  compact set $K\subset\Gw$, we assume that $\gs$ is a positive measure in $\Gw$ absolutely continuous with respect to the Lebesgue measure, with a nonnegative density $w$. For proving our results we will assume  that the function $\gw=w^{-\frac{1}{q-1}}$ is  $q'$-admissible in the sense of \cite[Chap 1]{HKM}. One sufficient condition is that $w$ belongs to the Muckenhoupt class $A_q$, that is 
\bel{E3-0}\BA {lll} \displaystyle
\sup_{B}\left(\frac{1}{\abs B}\myint{B}{}w dx\right)\left(\frac{1}{\abs B}\myint{B}{}w^{-\frac{1}{q-1}} dx\right)^{\frac{1}{p-1}}=m_{w,q}<\infty
\EA\ee
for all ball $B\subset \BBR^N$.\medskip

If $K\subset\Gw$ is compact, we set
\bel{E2}\BA {lll} 
c^\gw_{q}(K)=\inf\left\{\myint{\Gw}{}\abs{\Gd \gz}^{q'}\gw dx: \gz\in C^{\infty}_0(\Gw),\, \gz\geq 1\text{ in a neighborhood of }K\right\}.
\EA\ee
This defines a capacity on Borel subsets of $\Gw$. Since $\gw$ is $q'$-admissible, it satisfies Poincar\'e inequality, hence a set with zero $c^\gw_{q}$-capacity is $\gw$-negligible. Furthermore, following the proof of \cite[Th 3.3.3]{AdHe}, $c^\gw_{q}$ is equivalent to $\dot c^\gw_{q}$ defined by
\bel{E2-0}\BA {lll} 
\dot c^\gw_{q}(K)=\inf\left\{\norm{\gz}^{q'}_{W^{2,q'}_\gw}: \gz\in C^{\infty}_0(\Gw),\,0\leq\gz\leq 1,\,\gz\geq 1\text{ in a neighborhood of }K\right\}.
\EA\ee
The dual definition is ( see \cite[Th 2.5.1]{AdHe})
\bel{E2'}\BA {lll} 
\left(c^\gw_{q}(K)\right)^{\frac{1}{q'}}=\sup\left\{\gl(K): \gl\in \mathfrak M_b^+(K), \,\norm{\BBG[\gl]}_{L^q_\gw}\leq 1\right\}.
\EA\ee
\medskip

\nind{\it Proof of Theorem G}. {\it Step 1: The condition is sufficient.} We assume first that $L^q_{w,loc}(\Gw\setminus K)\cap u\in L^1(\Gw\setminus K)$ is a nonnegative subsolution of $(\ref{Z20})$ in the sense of distributions in $\Gw\setminus K$ where $K\subset\Gw$ is a compact subset with $c^\gw_{q}$-capacity zero. There exists a sequence of functions $\{\gz_k\}\subset C^{\infty}_0(\Gw)$ with value in $[0,1]$,  value $1$ in a neighborhood of $K$ and such that $\norm{\Gd\gz_k}_{L^{q'}_\gw}\to 0$ when $k\to\infty$. 
Let $\rho\in C^{\infty}_0(\Gw)$, $0\leq \rho\leq 1$, such that $\rho=1$ in a neighborhood of $K$ containing the support of the $\gz_k$. 
Using $\phi_k:=(1-\gz_k)^{\ga}\rho^{\ga}$, with $\ga>1$, in the very weak formulation of 
equation $(\ref{Z20})$ we obtain,  
\bel{E3}\BA {lll} \myint{\Gw}{}u^q \phi_k w dx\leq \myint{\Gw}{}u\Gd\phi_k dx\\[4mm]
\phantom{\myint{\Gw}{}u^q \phi_k w dx}
\leq\ga\myint{\Gw}{}u(1-\gz_k)^{\ga}\gr^{\ga-1}\Gd\rho dx-2\ga\myint{\Gw}{}u(1-\gz_k)^{\ga-1}\nabla\gz_k.\nabla\rho^{\ga} dx\\[4mm]
\phantom{\myint{\Gw}{}u^q \phi_k w dx}
-\alpha\myint{\Gw}{}u(1-\gz_k)^{\ga-1}\rho^{\alpha}\Gd\gz_k dx
+\ga(\ga-1)\myint{\Gw}{}u(1-\gz_k)^{\ga-2}\rho^{\ga}\abs{\nabla\gz_k}^2 dx\\[4mm]
\phantom{\myint{\Gw}{}u^q \phi_k w dx}
+\ga(\ga-1)\myint{\Gw}{}u(1-\gz_k)^{\ga}\rho^{\ga-2}\abs{\nabla\gr}^2 dx.
\EA\ee
Notice that  the second integral in the right-hand side vanishes since $\nabla\gz_k.\nabla\rho^{\ga}=0$ by the assumption on their support. 
If we choose $\ga=2q'$, we can bound the remaining integrals as follows: 
$$\BA {lll}\left|\myint{\Gw}{}u(1-\gz_k)^{2q'-1}\rho^{2q'}\Gd\gz_k dx\right|\leq \left(\myint{\Gw}{}u^q \phi_k w dx\right)^{\frac 1q}
\left(\myint{\Gw}{}\abs{\Gd \gz_k}^{q'}(1-\gz_k)^{q'}\gr^{2q'}\gw dx\right)^{\frac 1{q'}}\\[4mm]
\phantom{\left|\myint{\Gw}{}u(1-\gz_k)^{2q'-1}\rho^{2q'}\Gd\gz_k dx\right|}
\leq \left(\myint{\Gw}{}u^q \phi_k w dx\right)^{\frac 1q}
\left(\myint{\Gw}{}\abs{\Gd \gz_k}^{q'}\gw dx\right)^{\frac 1{q'}},
\EA$$
$$\BA {lll}
\left|\myint{\Gw}{}u(1-\gz_k)^{2q'}\gr^{2q'-1}\Gd\rho dx\right|
\leq \left(\myint{\Gw}{}u^q \phi_k w dx\right)^{\frac 1q}
\left(\myint{\Gw}{}\abs{\Gd\rho}^{q'} (1-\gz_k)^{2q'}\rho^{q'}\gw dx\right)^{\frac 1{q'}}\\[4mm]
\phantom{\left|\myint{\Gw}{}u(1-\gz_k)^{2q'}\gr^{2q'-1}\Gd\rho dx\right|}
\leq\left(\myint{\Gw}{}u^q \phi_k w dx\right)^{\frac 1q} 
\left(\myint{\Gw}{}\abs{\Gd\rho}^{q'} \gw dx\right)^{\frac 1{q'}},
\EA$$
$$\BA {lll}
\left|\myint{\Gw}{}u(1-\gz_k)^{2q'-2}\rho^{2q'}\abs{\nabla\gz_k}^2 dx\right|\leq  \left(\myint{\Gw}{}u^q \phi_k w dx\right)^{\frac 1q}
\left(\myint{\Gw}{}\abs{\nabla \gz_k}^{2q'}\gr^{2q'}\gw dx\right)^{\frac 1{q'}}\\[4mm]
\phantom{\left|\myint{\Gw}{}u(1-\gz_k)^{2q'-2}\rho^{2q'}\abs{\nabla\gz_k}^2 dx\right|}
\leq \left(\myint{\Gw}{}u^q \phi_k w dx\right)^{\frac 1q}
\left(\myint{\Gw}{}\abs{\nabla \gz_k}^{2q'}\gw dx\right)^{\frac 1{q'}},
\EA$$
and finally
$$\BA {lll}
\left|\myint{\Gw}{}u(1-\gz_k)^{2q'}\rho^{2q'-2}\abs{\nabla\rho}^2 dx\right|\leq  \left(\myint{\Gw}{}u^q \phi_k w dx\right)^{\frac 1q}
\left(\myint{\Gw}{}\abs{\nabla \rho}^{2q'}(1-\gz_k)^{2q'}\gw dx\right)^{\frac 1{q'}}\\[4mm]
\phantom{\left|\myint{\Gw}{}u(1-\gz_k)^{2q'-2}\rho^{2q'}\abs{\nabla\gz_k}^2 dx\right|}
\leq \left(\myint{\Gw}{}u^q \phi_k w dx\right)^{\frac 1q}
\left(\myint{\Gw}{}\abs{\nabla \rho}^{2q'}\gw dx\right)^{\frac 1{q'}}.
\EA$$
Because the Gagliardo-Nirenberg inequality holds with the $q'$-admissible weight $\gw$, we have for some $\gt\in (0,1)$ and some $c=c(q,N)>0$, 
\bel{E4}\BA {lll} 
\left(\myint{\Gw}{}\abs{\nabla \gz_k}^{2q'}\gw dx\right)^{\frac 1{2q'}}
\leq c\left(\myint{\Gw}{}\abs{\Gd \gz_k}^{q'}\gw dx\right)^{\frac \gt{q'}}\norm{\gz_k}_{L^\infty}^{1-\gt}
\\[4mm]
\phantom{\left(\myint{\Gw}{}\abs{\nabla \gz_k}^{2q'}\gw dx\right)^{\frac 1{q'}}}
\leq c'\left(\myint{\Gw}{}\abs{\Gd \gz_k}^{q'}\gw dx\right)^{\frac \gt{q'}}.
\EA\ee
Therefore, if we set 
$$X_k=\left(\myint{\Gw}{}u^q \phi_k w dx\right)^{\frac 1q}\quad\text{and }\;Z_k=\left(\myint{\Gw}{}\abs{\Gd \gz_k}^{q'}\gw dx\right)^{\frac 1{q'}},
$$
we obtain the inequation
\bel{E5}\BA {lll}
X_k^q\leq c_1X_kZ_k+c_2X_k+c_3X_kZ_k^\gt,
\EA\ee
for some positive constants $c_1,c_2,c_3$ depending on $q$, $N$ and $\gr$. 
By definition of $\zeta_k$ we have $Z_k\to 0$. We thus deduce that 
$X_k^q\leq cX_k$ with $q>1$ and then that the sequence $\{X_k\}$ is bounded. 
Since $\zeta_k\to 0$ almost everywhere, we have $\phi_k\to \rho^{2q'}$ almost everywhere. 
It then follows by Fatou's lemma that 
\bel{E6}\BA {lll}
\myint{\Gw}{}u^q \rho^{2q'}w dx\leq c. 
\EA\ee
We deduce that $u\in L^q_{w,loc}(\Gw)$. Since $\gw^{-\frac{q'}{q}}\in L^1_{loc}(\Gw)$, we obtain that $L^1_{loc}(\Gw)$ by H\"older's inequality.
If $u\in L^q_{w,loc}(\Gw\setminus K)\cap u\in L^1(\Gw\setminus K)$ is a distributional solution of $(\ref{Z20})$ in $\Gw\setminus K$, then 
 $|u|$ is a nonnegative subsolution with the same integrability constraints and we derive $u\in L^q_{w,loc}(\Gw)\cap  L^1_{loc}(\Gw)$.\\

If  $\gf\in C^\infty_0(\Gw)$, we take $\gf (1-\gz_k)^{2q'}$ for test function of equation $(\ref{Z20})$ in $\CD'(\Gw\setminus K)$,
$$ -\int_\Omega u \Delta(\phi(1-\zeta_k)^{2q'})\,dx 
+ \int_\Omega |u|^{q-1}u\phi(1-\zeta_k)^{2q'}w\,dx = 0. $$ 
Since $u\in L^q_{w,loc}(\Omega)$, $\phi$ has compact support, and $\zeta_k\to 0$ almost everywhere, we can pass to the limit as $k\to +\infty$ in the second integral using Lebesgue convergence theorem and obtain 
$$ \int_\Omega |u|^{q-1}u\phi(1-\zeta_k)^{2q'}w\,dx 
\to \int_\Omega |u|^{q-1}u\phi w\,dx. $$ 
Moreover  we can pass to the limit in the first integral expanding the laplacian. 
Using that $u\in L^1_{loc}(\Gw)$ and that 
$\Delta\zeta_k\to 0$ in $L^{q'}_\omega$, it is easy to prove from the previous computation that 
$$\myint{\Gw}{}u(1-\gz_k)^{q'}\Gd\phi dx\to \myint{\Gw}{}u\Gd\phi dx\quad\text{as }k\to\infty,
$$
and 
$$\displaystyle
\lim_{k\to\infty}
\myint{\Gw}{}u(1-\gz_k)^{2q'-1}\nabla\gz_k.\nabla\phi dx=0=\lim_{k\to\infty}\myint{\Gw}{}u(1-\gz_k)^{2q'-1}\phi\Gd\gz_k dx.
$$
Hence 
\bel{E7}\BA {lll}
- \myint{\Gw}{}u\Gd\phi dx+\myint{\Gw}{}u^q \gf  w dx=0
\EA\ee
\smallskip

\nind{\it Step 2: The condition is necessary.}
Let $K$ be a compact set with positive $c^\gw_{q}$-capacity. 
According to \cite{AdHe}[Th 2.5.3] there exists an extremal $\gm_k\in\mathfrak M_b^+(K)$ 
in the dual formulation $(\ref{E2'})$ of the capacity
According to Theorem D, problem $(\ref{D1})$ with $\gm=\gm_K$ admits a positive solution which is therefore a positive solution of $(\ref{E1})$.
\qeda

\medskip 

 \nind{\bf Aknowledgments} The authors have been supported by the MathAmsud program 13Math-03 QUESP with fundings from CNRS, Minist\`ere des Affaires \'Etrang\`eres et Europ\'eennes, CONICET and MINCyT.


\end{document}